\documentclass[review,dvipdfmx]{elsarticle}
\usepackage{lineno,hyperref}
\modulolinenumbers[5]

\journal{Journal of  Operations Management}
\newcommand{\makespan}{C^{\max}}
\usepackage{color}

\usepackage{amsmath}
\usepackage{amsfonts}
\newtheorem{thm}{Theorem}[section]
\newtheorem{cor}[thm]{Corollary}

\newcommand{\diff}{\mathrm{d}}

\newcommand{\shape}{expected terminate time }
\newcommand{\years}[1]{#1}
\newcommand{\Href}[1]{Fig.\ref{#1}}
\newcommand{\bea}{\begin{eqnarray}}
\newcommand{\eea}{\end{eqnarray}}
\newcommand{\eq}{&=&}
\newcommand{\nn}{\nonumber\\  }
\newcommand{\Sref}[1]{Eq.(\ref{#1})}
\renewcommand{\a}{\alpha}
\newcommand{\f}[2]{\frac{#1}{#2}}

\begin{document}

\begin{frontmatter}

\title{
Universality of Makespan\\
in Flowshop Scheduling Problem
}



\author[shinzatoaddress]{Takashi Shinzato\corref{mycorrespondingauthor}}
\cortext[mycorrespondingauthor]{Corresponding author}
\ead{takashi.shinzato@r.hit-u.ac.jp}

\author[kobayashiaddress]{Kei Kobayashi}
\ead{kei@math.keio.ac.jp}

\author[kakuaddress]{Ikou Kaku}
\ead{kakuikou@tcu.ac.jp}

\address[shinzatoaddress]{Mori Arinori Center for Higher Education and Global Mobility, Hitotsubashi University, Kunitachi, Tokyo, Japan}
\address[kobayashiaddress]{Faculty of Science and Technology, Keio 
University, Yokohama, Kanagawa, Japan}
\address[kakuaddress]{Graduate School of Environmental and Information Studies, Tokyo City 
University, Yokohama, Kanagawa, Japan}

\begin{abstract}
Makespan, which is defined as the time difference between the starting 
time and 
the terminate time of a sequence of jobs or tasks, as 
{the time to traverse}  a
belt conveyor system, is well known as one of the most 
important criteria in scheduling problems. It is often used by manufacturing firms in practice in order to improve the operational efficiency with respect to the order of job processing to be performed. 
It is known that the performance of a {machine} 
depends on the particular timing of the job processing 
even if the job processing order is fixed. That is, the performance of 
{a system with respect to} flowshop processing 
depends on the procedure of scheduling. In this present work, we first 
 discuss the relationship between makespan and several scheduling procedures in detail by using a small example and provide an algorithm 
 for deriving the makespan. Using our proposed algorithm, several 
 numerical experiments are examined so as to reveal the relationship 
 between the typical behavior of makespan and the position of the fiducial 
 machine, with respect to several distinguished  distributions of 
the processing time. We also discuss the behavior of makespan by using the 
 properties of the shape functions used  in the context of percolation theory. Our 
 contributions are firstly giving a detail discussion on the universality of 
 makespan in {flowshop problems}  and obtaining several novel properties of makespan, as follows: 
 (1) makespan possesses universality in the sense of being little affected by a change in the probability distribution of the 
 processing time, (2) makespan can be 
{decomposed into the sum of} two shape functions, and (3) makespan is less affected by the dispatching rule than by the scheduling procedure.
\end{abstract}


\end{frontmatter}

\section{Introduction}

Recently, in manufacturing, the importance of producing a wide 
variety of products in small quantities in order to meet the variety of 
customer needs is increasingly being recognized {\cite{
Malakooti
}. }
{Flowshop processing systems} are among the most widely used 
manufacturing systems for such production of a wide variety of products in 
small quantities. Furthermore, flowshop scheduling plays one of the most 
important roles in {product planning and,} furthermore, the corresponding scheduling problem is one of the most well-known 
scheduling problems in production as a special case of the jobshop scheduling problem. 
One significant feature  of flowshop scheduling is that the 
orders of job processing of whole machines {in a processing system }
are consistent, and it is known that the production capacity of {the processing system} can be maximally desterilized  with  little 
effort compared with the jobshop scheduling approach.

In the case that the number of machines {in the processing system} is 
two, Johnson's algorithm can be used to find the optimal scheduling of 
{the order that jobs are processed}
 in time $O(N\log N)$ where the number of jobs is $N$, and the optimality of the job processing order derived 
from Johnson's algorithm is mathematically guaranteed if the order of 
all jobs processed in the {processing system} is same as the order for the  flowshop scheduling problem \cite{Johnson}. 
This approach could be easily adapted to {other special situations}  by 
French \cite{French}. 
 However, it is known that the optimal order of 
the case of  three or more machines (except for certain specific cases) is 
difficult to solve, since the flowshop problem is NP-hard \cite{Garey}. 
Many researchers have used meta-heuristic 
algorithms and dispatching rules in order to overcome this difficulty of 
the scheduling problem \cite{Nowicki,Osman,Bierwirth}. For instance, 
Nowicki and Smutnicki proposed one of the most famous meta-heuristic 
algorithms, the tabu search algorithm,  for resolving the flowshop 
scheduling problem 
\cite{Nowicki}. 
Osman and Potts also investigated 
applying permutation  to flowshop scheduling problem 
using a simulated annealing algorithm and compare it with other meta-heuristic algorithms\cite{Osman}.
Moreover, Bierwirth and St$\ddot{\rm o}$ppler developed 
a genetic algorithm so as to 
analyze the flowshop scheduling problem \cite{Bierwirth}. However, 
because of the problem being an NP-hard problem \cite{Garey}, 
{the genetic algorithm was} not sufficient to solve completely the 
flowshop scheduling problem for all cases.

In recent decades, with respect to the optimality of a
{processing system} composed on more than three 
machines, the theory of constraints has drawn great attention. The theory of constraints is a 
management paradigm proposed by Goldratt \cite{Goldratt}. The mechanism of management scheduling  is explained in the form of a novel in his book. He also developed software called optimized production technology.
The core concepts of the theory of constraints are to 
improve the performance of the bottleneck machine in the {processing 
system} and to subordinate the other machines to the bottleneck 
machine. That is, 
in the theory of constraints,
it is necessary to implement an improvement of the performance of the bottleneck machine before anything else.
Improving the bottleneck machine 
as required by the theory of constraints is easy to do, but this does not always improve the total throughput, 
as described in detail below. Since there also {exist cases}  that 
the total throughput can be greatly improved 
{by improving the performance of a non-bottleneck machine}, 
from the unified viewpoint, we need to examine the asymptotical behavior of 
the production capacity of the {processing system} with respect to the 
starting point of scheduling of flowshop, both numerically and mathematically.

This paper is organized as follows; in the next section, we set up a flowshop scheduling model for simplicity of 
our discussion and explain briefly 
the relationship between the 
production capacity of the {processing system} 
and the machine location in the {processing system} which is the starting 
point of the scheduling with the help of a small example. Moreover, 
our motivation in this paper 
of deriving {the universal properties from some examples} is explained  and 
a flowshop scheduling algorithm is proposed in order to solve  
systematically {flowshop problems}. In section \ref{sec3}, 
numerical experiments show the validity of our proposed approaches 
under several different distribution models, namely, for a processing time 
independently and identically distributed according to an exponential, 
uniform, or $\chi^2$ distribution. We discuss mathematically the typical behaviors of the production 
capacity of the {processing system}  using the shape 
function developed in percolation theory in section \ref{sec4} and are able to explain the results of the numerical experiments mathematically. The last section is devoted to a summary and a statement of planned future work.

\begin{table}[t]
\begin{center}
\begin{tabular}{|c|c|c|}
\hline
Fiducial machine &processing time $y_\mu$ & makespan $\makespan$ \\
\hline
M1&$y_\mu=22$&$\makespan=48$\\
\hline
M2&$y_\mu=25$&$\makespan=49$\\
\hline
M3&$y_\mu=26$&$\makespan=50$\\
\hline
M4&$y_\mu=18$&$\makespan=55$\\
\hline
M5&$y_\mu=20$&$\makespan=48$\\
\hline
M6&$y_\mu=23$&$\makespan=48$\\
\hline
M7&$y_\mu=17$&$\makespan=48$\\
\hline
\end{tabular}
\caption{
\label{small-system-tab1}
The relationship between the makespan and  the position 
of the fiducial machine and processing time table is shown in table 
 \ref{small-system-tab2}, where the job processing order is fixed. The 
 Gantt charts are shown in 
\Href{gantt-chart1} 
to \Href{gantt-chart7}. 
}
\end{center}
\end{table}
\clearpage

\section{Model Setting and Flowshop Scheduling Algorithm}　

\subsection{Model Setting}
{We consider the flowshop scheduling problem of processing $N$ jobs on a
belt conveyor and $M$ machines configured in a line. 
$x_{\mu, 
i}$ indicates the processing time,} which includes the set-up time, of 
machine $\mu$ and job $i$, where $\mu=1,\cdots,M$ and $i=1,\cdots,N$, 
and the processing time depends on the workload of the particular job being processed and by the particular machine.
For convenience, the processing time is independently and identically distributed according to a given probability distribution ${\rm Pr}(x)$.  As the first step of our analytical research on flowshop 
scheduling, we assume that this {processing system} 
has no restriction on the due dates of the jobs. Further, 
$s_{\mu, i}$ and $t_{\mu, i}$ represent the starting time and terminate time of machine $\mu$ and job $i$, respectively, and satisfy the following relation:
\bea
t_{\mu,i}\eq s_{\mu,i}+x_{\mu,i}.\label{eq1}
\eea
Several different criteria for flowshop scheduling, for instance, average flow time, 
makespan, and average tardiness \cite{Malakooti}, have been discussed in 
the literature {\cite{Nawaz,Ho,Widmer,Hoos},} out of which 
the makespan under flowshop scheduling is focused on in the present paper. Here, makespan $\makespan$ means the total processing time of a {processing system} (e.g., a belt conveyor) and is defined as follows:
\bea
\makespan\eq t_{M,N}-s_{1,1}.
\eea
That is, the total processing time of the {processing system} is 
equal to the time difference between 
the starting time, $s_{1,1}$, and 
the terminate time, $t_{M,N}$, of the {processing system}. 

One of the goals of the flowshop scheduling problem is to minimize the makespan under 
permutations of the job processing order, 
{given the processing times of different jobs;}
 however, even if the processing order is  
fixed, as mentioned below, since the makespan depends on the position of 
the scheduling fiducial machine (where herein fiducial machine means the first machine in the processing scheduling) within the {processing system}. Namely, before the 
most discussed topic in the previous works (which is permuting job processing order) is analyzed, we need to 
reveal in detail the relationship between makespan and the position of the scheduling 
fiducial machine in the {processing system}. 
Although we start off handling only one 
example, in tables  \ref{small-system-tab1} and \ref{small-system-tab2}, based on the results, for a fixed job processing order, the 
processing scheduling by starting at the center machine of this {processing 
system} 
gives a lower processing capacity for this {processing system} than does that by starting at the bottleneck machine. 
We need to systematically discuss this property of makespan in order to 
analyze the optimal processing capacity before examining permutations of job processing order using the dispatching rules which are widely discussed in previous works.

\subsection{Example of Small System\label{section-of-example-of-small-system}}
We here discuss the relationship between makespan and the position of the scheduling fiducial machine 
in the {processing system} using the Gantt chart of a small example system. Let $y_\mu$ be the sum of the necessary processing times of machine $\mu$ if other machines do not influence machine $\mu$, that is, 
\bea
y_\mu\eq\sum_{i=1}^Nx_{\mu,i}.
\eea
 In accordance with common sense, 
the sum of {processing times in practice} of machine $\mu$ in the 
{processing system} is not smaller than the sum of necessary processing times $y_\mu$. From this, the bottleneck machine of this production system, 
as referred to in the literature on the theory of constraints,  is determined by the following:
\bea
\mu^*\eq\arg\mathop{\max}_{\mu}y_\mu.
\eea
Typically, the bottleneck machine is uniquely determined (namely, it is 
the machine which maximize the sum of the necessary 
processing times) if the processing time table is randomly {generated}.

As shown in table \ref{small-system-tab1}, M3 is the bottleneck machine in 
{this system}. At the first step, the processing scheduling of 
the bottleneck machine is independently determined using the processing time table.
At the next step, the processing 
scheduling of the nearest neighbors of the bottleneck machine, here 
M2({previous} machine in the {processing system}) and M4 ({next} machine), is 
uniquely assigned and depends only on the processing 
{schedules} of the bottleneck machine. In the same way, at the 
third step, the processing scheduling of M1, M5, M6, and M7 is decided in sequence. As a result, 
the Gantt chart of the bottleneck machine is as shown in {Fig.} \ref{gantt-chart3}. 
When the other machines are 
fiducial as the starting machine in scheduling, the Gantt charts are 
as shown in {Figs.} \ref{gantt-chart1}, 2, and 4 to 
\ref{gantt-chart7}. From these figures, we can see that the makespan of the bottleneck 
machine, here M3, is not smallest among the machines. Rather, the makespans of M1 and M7 (as well as M5 and M6), which are the first and last machines in 
this {processing system}, are minimal, whereas 
the makespan of M4, which is the center machine of {this system}, 
is maximal. This result conflicts with the canons of the theory of constraints, that is, 
the processing of the bottleneck machine
 is the constraint condition in this 
{processing
 system} and/or improvement (of the processing scheduling) of the 
bottleneck machine is most important. 
Our goal in this paper is to analyze the typical 
behaviors of makespan with respect to the position of the fiducial 
machine in the {processing system}.

\begin{table}[tb]
{\centering
\begin{tabular}{|c|c|c|c|c|c|c|c|c|c|c|}
\hline
&J1&J2&J3&J4&J5&J6&J7&J8&$y_\mu$&$\makespan$\\ \hline
M1&2&5&4&1&3&1&2&4&22&48\\ \hline
M2&2&3&2&4&3&4&2&5&25&49\\ \hline
M3&3&5&4&1&5&1&3&4&26&50\\ \hline
M4&4&5&1&1&3&1&2&1&18&55\\ \hline
M5&1&2&4&2&5&1&4&1&20&48\\ \hline
M6&5&2&3&1&4&1&5&2&23&48\\ \hline
M7&1&1&1&4&1&4&4&1&17&48\\ \hline
\end{tabular}
\caption{Processing time table of $N=8,M=7$. M3 is the bottleneck machine in this system and the total necessary processing time of the bottleneck machine 
is 26.  \label{small-system-tab2}}
}
\end{table}

\begin{figure}[bt]
{\centering
\begin{picture}(200,130)(-30,0)
\put(0,0){\vector(1,0){200}}
\put(0,0){\line(0,1){120}}
\put(-20,105){M1}
\put(-20,90){M2}
\put(-20,75){M3}
\put(-20,60){M4}
\put(-20,45){M5}
\put(-20,30){M6}
\put(-20,15){M7}
\textcolor{red}{\put(0,100){\dashbox{0.1}(6,15){\tiny 2}}\put(6,85){\dashbox{0.1}(6,15){\tiny 2}}\put(12,70){\dashbox{0.1}(9,15){\tiny 3}}\put(21,55){\dashbox{0.1}(12,15){\tiny 4}}\put(33,40){\dashbox{0.1}(3,15){\tiny 1}}\put(36,25){\dashbox{0.1}(15,15){\tiny 5}}\put(51,10){\dashbox{0.1}(3,15){\tiny 1}}}
\textcolor{blue}{\put(3,100){\dashbox{0.1}(15,15){\tiny 5}}\put(18,85){\dashbox{0.1}(9,15){\tiny 3}}\put(27,70){\dashbox{0.1}(15,15){\tiny 5}}\put(42,55){\dashbox{0.1}(15,15){\tiny 5}}\put(57,40){\dashbox{0.1}(6,15){\tiny 2}}\put(63,25){\dashbox{0.1}(6,15){\tiny 2}}\put(69,10){\dashbox{0.1}(3,15){\tiny 1}}}
\textcolor{magenta}{\put(15,100){\dashbox{0.1}(12,15){\tiny 4}}\put(27,85){\dashbox{0.1}(6,15){\tiny 2}}\put(39,70){\dashbox{0.1}(12,15){\tiny 4}}\put(54,55){\dashbox{0.1}(3,15){\tiny 1}}\put(60,40){\dashbox{0.1}(12,15){\tiny 4}}\put(72,25){\dashbox{0.1}(9,15){\tiny 3}}\put(81,10){\dashbox{0.1}(3,15){\tiny 1}}}
\textcolor{black}{\put(24,100){\dashbox{0.1}(3,15){\tiny 1}}\put(30,85){\dashbox{0.1}(12,15){\tiny 4}}\put(48,70){\dashbox{0.1}(3,15){\tiny 1}}\put(54,55){\dashbox{0.1}(3,15){\tiny 1}}\put(69,40){\dashbox{0.1}(6,15){\tiny 2}}\put(78,25){\dashbox{0.1}(3,15){\tiny 1}}\put(81,10){\dashbox{0.1}(12,15){\tiny 4}}}
\textcolor{red}{\put(24,100){\dashbox{0.1}(9,15){\tiny 3}}\put(39,85){\dashbox{0.1}(9,15){\tiny 3}}\put(48,70){\dashbox{0.1}(15,15){\tiny 5}}\put(63,55){\dashbox{0.1}(9,15){\tiny 3}}\put(72,40){\dashbox{0.1}(15,15){\tiny 5}}\put(87,25){\dashbox{0.1}(12,15){\tiny 4}}\put(99,10){\dashbox{0.1}(3,15){\tiny 1}}}
\textcolor{blue}{\put(30,100){\dashbox{0.1}(3,15){\tiny 1}}\put(45,85){\dashbox{0.1}(12,15){\tiny 4}}\put(60,70){\dashbox{0.1}(3,15){\tiny 1}}\put(69,55){\dashbox{0.1}(3,15){\tiny 1}}\put(84,40){\dashbox{0.1}(3,15){\tiny 1}}\put(96,25){\dashbox{0.1}(3,15){\tiny 1}}\put(99,10){\dashbox{0.1}(12,15){\tiny 4}}}
\textcolor{magenta}{\put(30,100){\dashbox{0.1}(6,15){\tiny 2}}\put(54,85){\dashbox{0.1}(6,15){\tiny 2}}\put(60,70){\dashbox{0.1}(9,15){\tiny 3}}\put(69,55){\dashbox{0.1}(6,15){\tiny 2}}\put(84,40){\dashbox{0.1}(12,15){\tiny 4}}\put(96,25){\dashbox{0.1}(15,15){\tiny 5}}\put(111,10){\dashbox{0.1}(12,15){\tiny 4}}}
\textcolor{black}{\put(33,100){\dashbox{0.1}(12,15){\tiny 4}}\put(57,85){\dashbox{0.1}(15,15){\tiny 5}}\put(72,70){\dashbox{0.1}(12,15){\tiny 4}}\put(84,55){\dashbox{0.1}(3,15){\tiny 1}}\put(93,40){\dashbox{0.1}(3,15){\tiny 1}}\put(108,25){\dashbox{0.1}(6,15){\tiny 2}}\put(120,10){\dashbox{0.1}(3,15){\tiny 1}}}
\end{picture}
\caption{{Gantt chart for M1 as the} fiducial machine. $\makespan=48$.\label{gantt-chart1}} 
}
{\centering
\begin{picture}(200,130)(-30,0)
\put(0,0){\vector(1,0){200}}
\put(0,0){\line(0,1){120}}
\put(-20,105){M1}
\put(-20,90){M2}
\put(-20,75){M3}
\put(-20,60){M4}
\put(-20,45){M5}
\put(-20,30){M6}
\put(-20,15){M7}
\textcolor{red}{\put(0,100){\dashbox{0.1}(6,15){\tiny 2}}\put(18,85){\dashbox{0.1}(6,15){\tiny 2}}\put(24,70){\dashbox{0.1}(9,15){\tiny 3}}\put(33,55){\dashbox{0.1}(12,15){\tiny 4}}\put(45,40){\dashbox{0.1}(3,15){\tiny 1}}\put(48,25){\dashbox{0.1}(15,15){\tiny 5}}\put(63,10){\dashbox{0.1}(3,15){\tiny 1}}}
\textcolor{blue}{\put(3,100){\dashbox{0.1}(15,15){\tiny 5}}\put(21,85){\dashbox{0.1}(9,15){\tiny 3}}\put(30,70){\dashbox{0.1}(15,15){\tiny 5}}\put(45,55){\dashbox{0.1}(15,15){\tiny 5}}\put(60,40){\dashbox{0.1}(6,15){\tiny 2}}\put(66,25){\dashbox{0.1}(6,15){\tiny 2}}\put(72,10){\dashbox{0.1}(3,15){\tiny 1}}}
\textcolor{magenta}{\put(15,100){\dashbox{0.1}(12,15){\tiny 4}}\put(27,85){\dashbox{0.1}(6,15){\tiny 2}}\put(42,70){\dashbox{0.1}(12,15){\tiny 4}}\put(57,55){\dashbox{0.1}(3,15){\tiny 1}}\put(63,40){\dashbox{0.1}(12,15){\tiny 4}}\put(75,25){\dashbox{0.1}(9,15){\tiny 3}}\put(84,10){\dashbox{0.1}(3,15){\tiny 1}}}
\textcolor{black}{\put(27,100){\dashbox{0.1}(3,15){\tiny 1}}\put(30,85){\dashbox{0.1}(12,15){\tiny 4}}\put(51,70){\dashbox{0.1}(3,15){\tiny 1}}\put(57,55){\dashbox{0.1}(3,15){\tiny 1}}\put(72,40){\dashbox{0.1}(6,15){\tiny 2}}\put(81,25){\dashbox{0.1}(3,15){\tiny 1}}\put(84,10){\dashbox{0.1}(12,15){\tiny 4}}}
\textcolor{red}{\put(30,100){\dashbox{0.1}(9,15){\tiny 3}}\put(39,85){\dashbox{0.1}(9,15){\tiny 3}}\put(51,70){\dashbox{0.1}(15,15){\tiny 5}}\put(66,55){\dashbox{0.1}(9,15){\tiny 3}}\put(75,40){\dashbox{0.1}(15,15){\tiny 5}}\put(90,25){\dashbox{0.1}(12,15){\tiny 4}}\put(102,10){\dashbox{0.1}(3,15){\tiny 1}}}
\textcolor{blue}{\put(42,100){\dashbox{0.1}(3,15){\tiny 1}}\put(45,85){\dashbox{0.1}(12,15){\tiny 4}}\put(63,70){\dashbox{0.1}(3,15){\tiny 1}}\put(72,55){\dashbox{0.1}(3,15){\tiny 1}}\put(87,40){\dashbox{0.1}(3,15){\tiny 1}}\put(99,25){\dashbox{0.1}(3,15){\tiny 1}}\put(102,10){\dashbox{0.1}(12,15){\tiny 4}}}
\textcolor{magenta}{\put(42,100){\dashbox{0.1}(6,15){\tiny 2}}\put(54,85){\dashbox{0.1}(6,15){\tiny 2}}\put(63,70){\dashbox{0.1}(9,15){\tiny 3}}\put(72,55){\dashbox{0.1}(6,15){\tiny 2}}\put(87,40){\dashbox{0.1}(12,15){\tiny 4}}\put(99,25){\dashbox{0.1}(15,15){\tiny 5}}\put(114,10){\dashbox{0.1}(12,15){\tiny 4}}}
\textcolor{black}{\put(45,100){\dashbox{0.1}(12,15){\tiny 4}}\put(57,85){\dashbox{0.1}(15,15){\tiny 5}}\put(72,70){\dashbox{0.1}(12,15){\tiny 4}}\put(84,55){\dashbox{0.1}(3,15){\tiny 1}}\put(96,40){\dashbox{0.1}(3,15){\tiny 1}}\put(111,25){\dashbox{0.1}(6,15){\tiny 2}}\put(123,10){\dashbox{0.1}(3,15){\tiny 1}}}
\end{picture}
\caption{{Gantt chart for M2 as the} fiducial machine. $\makespan=49$.\label{gantt-chart2}}
}
\end{figure}

\begin{figure}[bt]
{\centering
\begin{picture}(200,130)(-30,0)
\put(0,0){\vector(1,0){200}}
\put(0,0){\line(0,1){120}}
\put(-20,105){M1}
\put(-20,90){M2}
\put(-20,75){M3}
\put(-20,60){M4}
\put(-20,45){M5}
\put(-20,30){M6}
\put(-20,15){M7}
\textcolor{red}{\put(0,100){\dashbox{0.1}(6,15){\tiny 2}}\put(18,85){\dashbox{0.1}(6,15){\tiny 2}}\put(27,70){\dashbox{0.1}(9,15){\tiny 3}}\put(36,55){\dashbox{0.1}(12,15){\tiny 4}}\put(48,40){\dashbox{0.1}(3,15){\tiny 1}}\put(51,25){\dashbox{0.1}(15,15){\tiny 5}}\put(66,10){\dashbox{0.1}(3,15){\tiny 1}}}
\textcolor{blue}{\put(3,100){\dashbox{0.1}(15,15){\tiny 5}}\put(21,85){\dashbox{0.1}(9,15){\tiny 3}}\put(33,70){\dashbox{0.1}(15,15){\tiny 5}}\put(48,55){\dashbox{0.1}(15,15){\tiny 5}}\put(63,40){\dashbox{0.1}(6,15){\tiny 2}}\put(69,25){\dashbox{0.1}(6,15){\tiny 2}}\put(75,10){\dashbox{0.1}(3,15){\tiny 1}}}
\textcolor{magenta}{\put(15,100){\dashbox{0.1}(12,15){\tiny 4}}\put(27,85){\dashbox{0.1}(6,15){\tiny 2}}\put(45,70){\dashbox{0.1}(12,15){\tiny 4}}\put(60,55){\dashbox{0.1}(3,15){\tiny 1}}\put(66,40){\dashbox{0.1}(12,15){\tiny 4}}\put(78,25){\dashbox{0.1}(9,15){\tiny 3}}\put(87,10){\dashbox{0.1}(3,15){\tiny 1}}}
\textcolor{black}{\put(27,100){\dashbox{0.1}(3,15){\tiny 1}}\put(30,85){\dashbox{0.1}(12,15){\tiny 4}}\put(54,70){\dashbox{0.1}(3,15){\tiny 1}}\put(60,55){\dashbox{0.1}(3,15){\tiny 1}}\put(75,40){\dashbox{0.1}(6,15){\tiny 2}}\put(84,25){\dashbox{0.1}(3,15){\tiny 1}}\put(87,10){\dashbox{0.1}(12,15){\tiny 4}}}
\textcolor{red}{\put(30,100){\dashbox{0.1}(9,15){\tiny 3}}\put(39,85){\dashbox{0.1}(9,15){\tiny 3}}\put(54,70){\dashbox{0.1}(15,15){\tiny 5}}\put(69,55){\dashbox{0.1}(9,15){\tiny 3}}\put(78,40){\dashbox{0.1}(15,15){\tiny 5}}\put(93,25){\dashbox{0.1}(12,15){\tiny 4}}\put(105,10){\dashbox{0.1}(3,15){\tiny 1}}}
\textcolor{blue}{\put(42,100){\dashbox{0.1}(3,15){\tiny 1}}\put(45,85){\dashbox{0.1}(12,15){\tiny 4}}\put(66,70){\dashbox{0.1}(3,15){\tiny 1}}\put(75,55){\dashbox{0.1}(3,15){\tiny 1}}\put(90,40){\dashbox{0.1}(3,15){\tiny 1}}\put(102,25){\dashbox{0.1}(3,15){\tiny 1}}\put(105,10){\dashbox{0.1}(12,15){\tiny 4}}}
\textcolor{magenta}{\put(42,100){\dashbox{0.1}(6,15){\tiny 2}}\put(54,85){\dashbox{0.1}(6,15){\tiny 2}}\put(66,70){\dashbox{0.1}(9,15){\tiny 3}}\put(75,55){\dashbox{0.1}(6,15){\tiny 2}}\put(90,40){\dashbox{0.1}(12,15){\tiny 4}}\put(102,25){\dashbox{0.1}(15,15){\tiny 5}}\put(117,10){\dashbox{0.1}(12,15){\tiny 4}}}
\textcolor{black}{\put(45,100){\dashbox{0.1}(12,15){\tiny 4}}\put(57,85){\dashbox{0.1}(15,15){\tiny 5}}\put(72,70){\dashbox{0.1}(12,15){\tiny 4}}\put(84,55){\dashbox{0.1}(3,15){\tiny 1}}\put(99,40){\dashbox{0.1}(3,15){\tiny 1}}\put(114,25){\dashbox{0.1}(6,15){\tiny 2}}\put(126,10){\dashbox{0.1}(3,15){\tiny 1}}}
\end{picture}
\caption{{Gantt chart for M3 as the} fiducial machine. $\makespan=50$.\label{gantt-chart3}}
}
{\centering
\begin{picture}(200,130)(-30,0)
\put(0,0){\vector(1,0){200}}
\put(0,0){\line(0,1){120}}
\put(-20,105){M1}
\put(-20,90){M2}
\put(-20,75){M3}
\put(-20,60){M4}
\put(-20,45){M5}
\put(-20,30){M6}
\put(-20,15){M7}
\textcolor{red}{\put(0,100){\dashbox{0.1}(6,15){\tiny 2}}\put(18,85){\dashbox{0.1}(6,15){\tiny 2}}\put(27,70){\dashbox{0.1}(9,15){\tiny 3}}\put(54,55){\dashbox{0.1}(12,15){\tiny 4}}\put(66,40){\dashbox{0.1}(3,15){\tiny 1}}\put(69,25){\dashbox{0.1}(15,15){\tiny 5}}\put(84,10){\dashbox{0.1}(3,15){\tiny 1}}}
\textcolor{blue}{\put(3,100){\dashbox{0.1}(15,15){\tiny 5}}\put(21,85){\dashbox{0.1}(9,15){\tiny 3}}\put(33,70){\dashbox{0.1}(15,15){\tiny 5}}\put(63,55){\dashbox{0.1}(15,15){\tiny 5}}\put(78,40){\dashbox{0.1}(6,15){\tiny 2}}\put(84,25){\dashbox{0.1}(6,15){\tiny 2}}\put(90,10){\dashbox{0.1}(3,15){\tiny 1}}}
\textcolor{magenta}{\put(15,100){\dashbox{0.1}(12,15){\tiny 4}}\put(27,85){\dashbox{0.1}(6,15){\tiny 2}}\put(45,70){\dashbox{0.1}(12,15){\tiny 4}}\put(75,55){\dashbox{0.1}(3,15){\tiny 1}}\put(81,40){\dashbox{0.1}(12,15){\tiny 4}}\put(93,25){\dashbox{0.1}(9,15){\tiny 3}}\put(102,10){\dashbox{0.1}(3,15){\tiny 1}}}
\textcolor{black}{\put(27,100){\dashbox{0.1}(3,15){\tiny 1}}\put(30,85){\dashbox{0.1}(12,15){\tiny 4}}\put(54,70){\dashbox{0.1}(3,15){\tiny 1}}\put(75,55){\dashbox{0.1}(3,15){\tiny 1}}\put(90,40){\dashbox{0.1}(6,15){\tiny 2}}\put(99,25){\dashbox{0.1}(3,15){\tiny 1}}\put(102,10){\dashbox{0.1}(12,15){\tiny 4}}}
\textcolor{red}{\put(30,100){\dashbox{0.1}(9,15){\tiny 3}}\put(39,85){\dashbox{0.1}(9,15){\tiny 3}}\put(54,70){\dashbox{0.1}(15,15){\tiny 5}}\put(75,55){\dashbox{0.1}(9,15){\tiny 3}}\put(93,40){\dashbox{0.1}(15,15){\tiny 5}}\put(108,25){\dashbox{0.1}(12,15){\tiny 4}}\put(120,10){\dashbox{0.1}(3,15){\tiny 1}}}
\textcolor{blue}{\put(42,100){\dashbox{0.1}(3,15){\tiny 1}}\put(45,85){\dashbox{0.1}(12,15){\tiny 4}}\put(66,70){\dashbox{0.1}(3,15){\tiny 1}}\put(81,55){\dashbox{0.1}(3,15){\tiny 1}}\put(105,40){\dashbox{0.1}(3,15){\tiny 1}}\put(117,25){\dashbox{0.1}(3,15){\tiny 1}}\put(120,10){\dashbox{0.1}(12,15){\tiny 4}}}
\textcolor{magenta}{\put(42,100){\dashbox{0.1}(6,15){\tiny 2}}\put(54,85){\dashbox{0.1}(6,15){\tiny 2}}\put(66,70){\dashbox{0.1}(9,15){\tiny 3}}\put(81,55){\dashbox{0.1}(6,15){\tiny 2}}\put(105,40){\dashbox{0.1}(12,15){\tiny 4}}\put(117,25){\dashbox{0.1}(15,15){\tiny 5}}\put(132,10){\dashbox{0.1}(12,15){\tiny 4}}}
\textcolor{black}{\put(45,100){\dashbox{0.1}(12,15){\tiny 4}}\put(57,85){\dashbox{0.1}(15,15){\tiny 5}}\put(72,70){\dashbox{0.1}(12,15){\tiny 4}}\put(84,55){\dashbox{0.1}(3,15){\tiny 1}}\put(114,40){\dashbox{0.1}(3,15){\tiny 1}}\put(129,25){\dashbox{0.1}(6,15){\tiny 2}}\put(141,10){\dashbox{0.1}(3,15){\tiny 1}}}
\end{picture}
\caption{{Gantt chart for M4 as the} fiducial machine. $\makespan=55$.\label{gantt-chart4}}
}
\end{figure}

\begin{figure}[t]
{\centering
\begin{picture}(200,130)(-30,0)
\put(0,0){\vector(1,0){200}}
\put(0,0){\line(0,1){120}}
\put(-20,105){M1}
\put(-20,90){M2}
\put(-20,75){M3}
\put(-20,60){M4}
\put(-20,45){M5}
\put(-20,30){M6}
\put(-20,15){M7}
\textcolor{red}{\put(0,100){\dashbox{0.1}(6,15){\tiny 2}}\put(15,85){\dashbox{0.1}(6,15){\tiny 2}}\put(21,70){\dashbox{0.1}(9,15){\tiny 3}}\put(33,55){\dashbox{0.1}(12,15){\tiny 4}}\put(57,40){\dashbox{0.1}(3,15){\tiny 1}}\put(60,25){\dashbox{0.1}(15,15){\tiny 5}}\put(75,10){\dashbox{0.1}(3,15){\tiny 1}}}
\textcolor{blue}{\put(3,100){\dashbox{0.1}(15,15){\tiny 5}}\put(18,85){\dashbox{0.1}(9,15){\tiny 3}}\put(27,70){\dashbox{0.1}(15,15){\tiny 5}}\put(42,55){\dashbox{0.1}(15,15){\tiny 5}}\put(57,40){\dashbox{0.1}(6,15){\tiny 2}}\put(72,25){\dashbox{0.1}(6,15){\tiny 2}}\put(78,10){\dashbox{0.1}(3,15){\tiny 1}}}
\textcolor{magenta}{\put(15,100){\dashbox{0.1}(12,15){\tiny 4}}\put(27,85){\dashbox{0.1}(6,15){\tiny 2}}\put(39,70){\dashbox{0.1}(12,15){\tiny 4}}\put(57,55){\dashbox{0.1}(3,15){\tiny 1}}\put(60,40){\dashbox{0.1}(12,15){\tiny 4}}\put(75,25){\dashbox{0.1}(9,15){\tiny 3}}\put(84,10){\dashbox{0.1}(3,15){\tiny 1}}}
\textcolor{black}{\put(27,100){\dashbox{0.1}(3,15){\tiny 1}}\put(30,85){\dashbox{0.1}(12,15){\tiny 4}}\put(48,70){\dashbox{0.1}(3,15){\tiny 1}}\put(63,55){\dashbox{0.1}(3,15){\tiny 1}}\put(69,40){\dashbox{0.1}(6,15){\tiny 2}}\put(81,25){\dashbox{0.1}(3,15){\tiny 1}}\put(84,10){\dashbox{0.1}(12,15){\tiny 4}}}
\textcolor{red}{\put(30,100){\dashbox{0.1}(9,15){\tiny 3}}\put(39,85){\dashbox{0.1}(9,15){\tiny 3}}\put(48,70){\dashbox{0.1}(15,15){\tiny 5}}\put(63,55){\dashbox{0.1}(9,15){\tiny 3}}\put(72,40){\dashbox{0.1}(15,15){\tiny 5}}\put(87,25){\dashbox{0.1}(12,15){\tiny 4}}\put(99,10){\dashbox{0.1}(3,15){\tiny 1}}}
\textcolor{blue}{\put(48,100){\dashbox{0.1}(3,15){\tiny 1}}\put(51,85){\dashbox{0.1}(12,15){\tiny 4}}\put(69,70){\dashbox{0.1}(3,15){\tiny 1}}\put(78,55){\dashbox{0.1}(3,15){\tiny 1}}\put(84,40){\dashbox{0.1}(3,15){\tiny 1}}\put(96,25){\dashbox{0.1}(3,15){\tiny 1}}\put(99,10){\dashbox{0.1}(12,15){\tiny 4}}}
\textcolor{magenta}{\put(48,100){\dashbox{0.1}(6,15){\tiny 2}}\put(60,85){\dashbox{0.1}(6,15){\tiny 2}}\put(69,70){\dashbox{0.1}(9,15){\tiny 3}}\put(78,55){\dashbox{0.1}(6,15){\tiny 2}}\put(84,40){\dashbox{0.1}(12,15){\tiny 4}}\put(96,25){\dashbox{0.1}(15,15){\tiny 5}}\put(111,10){\dashbox{0.1}(12,15){\tiny 4}}}
\textcolor{black}{\put(51,100){\dashbox{0.1}(12,15){\tiny 4}}\put(63,85){\dashbox{0.1}(15,15){\tiny 5}}\put(78,70){\dashbox{0.1}(12,15){\tiny 4}}\put(90,55){\dashbox{0.1}(3,15){\tiny 1}}\put(93,40){\dashbox{0.1}(3,15){\tiny 1}}\put(108,25){\dashbox{0.1}(6,15){\tiny 2}}\put(120,10){\dashbox{0.1}(3,15){\tiny 1}}}
\end{picture}
\caption{{Gantt chart for M5 as the} fiducial machine. $\makespan=48$.\label{gantt-chart5}}
}
{\centering
\begin{picture}(200,130)(-30,0)
\put(0,0){\vector(1,0){200}}
\put(0,0){\line(0,1){120}}
\put(-20,105){M1}
\put(-20,90){M2}
\put(-20,75){M3}
\put(-20,60){M4}
\put(-20,45){M5}
\put(-20,30){M6}
\put(-20,15){M7}
\textcolor{red}{\put(0,100){\dashbox{0.1}(6,15){\tiny 2}}\put(15,85){\dashbox{0.1}(6,15){\tiny 2}}\put(21,70){\dashbox{0.1}(9,15){\tiny 3}}\put(33,55){\dashbox{0.1}(12,15){\tiny 4}}\put(57,40){\dashbox{0.1}(3,15){\tiny 1}}\put(66,25){\dashbox{0.1}(15,15){\tiny 5}}\put(81,10){\dashbox{0.1}(3,15){\tiny 1}}}
\textcolor{blue}{\put(3,100){\dashbox{0.1}(15,15){\tiny 5}}\put(18,85){\dashbox{0.1}(9,15){\tiny 3}}\put(27,70){\dashbox{0.1}(15,15){\tiny 5}}\put(42,55){\dashbox{0.1}(15,15){\tiny 5}}\put(57,40){\dashbox{0.1}(6,15){\tiny 2}}\put(78,25){\dashbox{0.1}(6,15){\tiny 2}}\put(84,10){\dashbox{0.1}(3,15){\tiny 1}}}
\textcolor{magenta}{\put(15,100){\dashbox{0.1}(12,15){\tiny 4}}\put(27,85){\dashbox{0.1}(6,15){\tiny 2}}\put(39,70){\dashbox{0.1}(12,15){\tiny 4}}\put(57,55){\dashbox{0.1}(3,15){\tiny 1}}\put(60,40){\dashbox{0.1}(12,15){\tiny 4}}\put(81,25){\dashbox{0.1}(9,15){\tiny 3}}\put(90,10){\dashbox{0.1}(3,15){\tiny 1}}}
\textcolor{black}{\put(27,100){\dashbox{0.1}(3,15){\tiny 1}}\put(30,85){\dashbox{0.1}(12,15){\tiny 4}}\put(48,70){\dashbox{0.1}(3,15){\tiny 1}}\put(63,55){\dashbox{0.1}(3,15){\tiny 1}}\put(69,40){\dashbox{0.1}(6,15){\tiny 2}}\put(87,25){\dashbox{0.1}(3,15){\tiny 1}}\put(90,10){\dashbox{0.1}(12,15){\tiny 4}}}
\textcolor{red}{\put(30,100){\dashbox{0.1}(9,15){\tiny 3}}\put(39,85){\dashbox{0.1}(9,15){\tiny 3}}\put(48,70){\dashbox{0.1}(15,15){\tiny 5}}\put(63,55){\dashbox{0.1}(9,15){\tiny 3}}\put(72,40){\dashbox{0.1}(15,15){\tiny 5}}\put(87,25){\dashbox{0.1}(12,15){\tiny 4}}\put(99,10){\dashbox{0.1}(3,15){\tiny 1}}}
\textcolor{blue}{\put(51,100){\dashbox{0.1}(3,15){\tiny 1}}\put(54,85){\dashbox{0.1}(12,15){\tiny 4}}\put(69,70){\dashbox{0.1}(3,15){\tiny 1}}\put(78,55){\dashbox{0.1}(3,15){\tiny 1}}\put(84,40){\dashbox{0.1}(3,15){\tiny 1}}\put(96,25){\dashbox{0.1}(3,15){\tiny 1}}\put(99,10){\dashbox{0.1}(12,15){\tiny 4}}}
\textcolor{magenta}{\put(57,100){\dashbox{0.1}(6,15){\tiny 2}}\put(63,85){\dashbox{0.1}(6,15){\tiny 2}}\put(69,70){\dashbox{0.1}(9,15){\tiny 3}}\put(78,55){\dashbox{0.1}(6,15){\tiny 2}}\put(84,40){\dashbox{0.1}(12,15){\tiny 4}}\put(96,25){\dashbox{0.1}(15,15){\tiny 5}}\put(111,10){\dashbox{0.1}(12,15){\tiny 4}}}
\textcolor{black}{\put(63,100){\dashbox{0.1}(12,15){\tiny 4}}\put(75,85){\dashbox{0.1}(15,15){\tiny 5}}\put(90,70){\dashbox{0.1}(12,15){\tiny 4}}\put(102,55){\dashbox{0.1}(3,15){\tiny 1}}\put(105,40){\dashbox{0.1}(3,15){\tiny 1}}\put(108,25){\dashbox{0.1}(6,15){\tiny 2}}\put(120,10){\dashbox{0.1}(3,15){\tiny 1}}}
\end{picture}
\caption{{Gantt chart for M6 as the} fiducial machine. $\makespan=48$.\label{gantt-chart6}}
}
{\centering
\begin{picture}(200,130)(-30,0)
\put(0,0){\vector(1,0){200}}
\put(0,0){\line(0,1){120}}
\put(-20,105){M1}
\put(-20,90){M2}
\put(-20,75){M3}
\put(-20,60){M4}
\put(-20,45){M5}
\put(-20,30){M6}
\put(-20,15){M7}
\textcolor{red}{\put(0,100){\dashbox{0.1}(6,15){\tiny 2}}\put(15,85){\dashbox{0.1}(6,15){\tiny 2}}\put(21,70){\dashbox{0.1}(9,15){\tiny 3}}\put(33,55){\dashbox{0.1}(12,15){\tiny 4}}\put(57,40){\dashbox{0.1}(3,15){\tiny 1}}\put(66,25){\dashbox{0.1}(15,15){\tiny 5}}\put(93,10){\dashbox{0.1}(3,15){\tiny 1}}}
\textcolor{blue}{\put(3,100){\dashbox{0.1}(15,15){\tiny 5}}\put(18,85){\dashbox{0.1}(9,15){\tiny 3}}\put(27,70){\dashbox{0.1}(15,15){\tiny 5}}\put(42,55){\dashbox{0.1}(15,15){\tiny 5}}\put(57,40){\dashbox{0.1}(6,15){\tiny 2}}\put(78,25){\dashbox{0.1}(6,15){\tiny 2}}\put(93,10){\dashbox{0.1}(3,15){\tiny 1}}}
\textcolor{magenta}{\put(15,100){\dashbox{0.1}(12,15){\tiny 4}}\put(27,85){\dashbox{0.1}(6,15){\tiny 2}}\put(39,70){\dashbox{0.1}(12,15){\tiny 4}}\put(57,55){\dashbox{0.1}(3,15){\tiny 1}}\put(60,40){\dashbox{0.1}(12,15){\tiny 4}}\put(81,25){\dashbox{0.1}(9,15){\tiny 3}}\put(93,10){\dashbox{0.1}(3,15){\tiny 1}}}
\textcolor{black}{\put(27,100){\dashbox{0.1}(3,15){\tiny 1}}\put(30,85){\dashbox{0.1}(12,15){\tiny 4}}\put(48,70){\dashbox{0.1}(3,15){\tiny 1}}\put(63,55){\dashbox{0.1}(3,15){\tiny 1}}\put(69,40){\dashbox{0.1}(6,15){\tiny 2}}\put(87,25){\dashbox{0.1}(3,15){\tiny 1}}\put(93,10){\dashbox{0.1}(12,15){\tiny 4}}}
\textcolor{red}{\put(30,100){\dashbox{0.1}(9,15){\tiny 3}}\put(39,85){\dashbox{0.1}(9,15){\tiny 3}}\put(48,70){\dashbox{0.1}(15,15){\tiny 5}}\put(63,55){\dashbox{0.1}(9,15){\tiny 3}}\put(72,40){\dashbox{0.1}(15,15){\tiny 5}}\put(87,25){\dashbox{0.1}(12,15){\tiny 4}}\put(102,10){\dashbox{0.1}(3,15){\tiny 1}}}
\textcolor{blue}{\put(51,100){\dashbox{0.1}(3,15){\tiny 1}}\put(54,85){\dashbox{0.1}(12,15){\tiny 4}}\put(69,70){\dashbox{0.1}(3,15){\tiny 1}}\put(78,55){\dashbox{0.1}(3,15){\tiny 1}}\put(84,40){\dashbox{0.1}(3,15){\tiny 1}}\put(96,25){\dashbox{0.1}(3,15){\tiny 1}}\put(102,10){\dashbox{0.1}(12,15){\tiny 4}}}
\textcolor{magenta}{\put(57,100){\dashbox{0.1}(6,15){\tiny 2}}\put(63,85){\dashbox{0.1}(6,15){\tiny 2}}\put(69,70){\dashbox{0.1}(9,15){\tiny 3}}\put(78,55){\dashbox{0.1}(6,15){\tiny 2}}\put(84,40){\dashbox{0.1}(12,15){\tiny 4}}\put(96,25){\dashbox{0.1}(15,15){\tiny 5}}\put(111,10){\dashbox{0.1}(12,15){\tiny 4}}}
\textcolor{black}{\put(69,100){\dashbox{0.1}(12,15){\tiny 4}}\put(81,85){\dashbox{0.1}(15,15){\tiny 5}}\put(96,70){\dashbox{0.1}(12,15){\tiny 4}}\put(108,55){\dashbox{0.1}(3,15){\tiny 1}}\put(111,40){\dashbox{0.1}(3,15){\tiny 1}}\put(114,25){\dashbox{0.1}(6,15){\tiny 2}}\put(120,10){\dashbox{0.1}(3,15){\tiny 1}}}
\end{picture}
\caption{{Gantt chart for M7 as the} fiducial machine. $\makespan=48$.\label{gantt-chart7}}
}
\end{figure}

\clearpage

\subsection{Forward Scheduling}
In the previous subsection, 
makespan of each machine in the {processing system} was discussed for an example of a small 
{processing system}. {Here,} the forward scheduling procedure is introduced in detail.

With respect to machine $\mu$, after the starting and terminate times of machine $\mu-1$ have been determined, the starting and terminate times of machine $\mu$ and job $i$, $s_{\mu,i}$ and $t_{\mu,i}$, are assigned as follows:
\bea
s_{\mu,i}\eq\max\left(
t_{\mu,i-1},
t_{\mu-1,i}
\right),\label{eq5}\\
t_{\mu,i}\eq s_{\mu,i}+x_{\mu,i},\label{eq6}
\eea
where the starting time $s_{\mu,i}$ is decided as the larger of 
the terminate time of machine $\mu$ and job $i-1$ and that of machine 
$\mu-1$ and job $i$, and the terminate time $t_{\mu,i}$ is assessed using 
equation (\ref{eq1}) or (\ref{eq6}).

\subsection{Backward Scheduling}
In a similar way to those in previous subsections, the backward scheduling procedure 
is explained in detail as follows. With respect to machine $\mu$, after those of machine $\mu+1$ have been assigned, the terminate and starting times of machine $\mu$ and job $i$, $t_{\mu,i}$ and $s_{\mu,i}$, are calculated as follows:
\bea
\label{eq7}
t_{\mu,i}\eq\min\left(s_{\mu,i+1},s_{\mu+1,i}\right),\\
s_{\mu,i}\eq t_{\mu,i}-x_{\mu,i},
\label{eq8}
\eea
where the terminate time $t_{\mu,i}$ is evaluated as the smaller of the starting time of 
machine $\mu$ and job $i+1$ and that of machine $\mu+1$ and job $i$, and 
the starting time $s_{\mu,i}$ is assessed using equation (\ref{eq1}) 
or (\ref{eq8}).

\subsection{Algorithm for Evaluating Makespan\label{sec2.5}}

Summarizing our {explanations} in the previous subsections, the algorithm for 
evaluating makespan with each choice of fiducial machine 
(the starting machine in terms of scheduling) in the {processing system} is organized as follows:
\begin{description}
\item[Step 0:]
The processing time of machine $\mu$ and job $i$, $x_{\mu,i}$, is 
	   randomly assigned according to probability distribution ${\rm Pr}(x)$ to create a processing 
	   time table matrix $X=\left\{x_{\mu,i}\right\}\in\mathbb{R}^{M\times N}$. Initially, the number of the fiducial machine is $\nu=1$.
\item[Step 1:]
With respect to fiducial machine $\nu$, the starting time of job $1$ is  
	   $s_{\nu,1}(\nu)=0$ by assumption. Moreover, since 
the starting and terminate times of the fiducial machine do not depend on the processing operation schedules of the other machines, the starting and terminate times of fiducial machine $\nu$ and job $i$ are determined in accordance with the following relations:
\bea
t_{\nu,i}(\nu)\eq s_{\nu,i}(\nu)+x_{\nu,i},\\
s_{\nu,i+1}(\nu)\eq t_{\nu,i}(\nu).
\eea
In Step 1, 	   
$s_{\nu,1}(\nu),t_{\nu,1}(\nu),\cdots,s_{\nu,N}(\nu),t_{\nu,N}(\nu)$ are 
	   uniquely {determined}. Note that hereafter we 
use $s_{\mu,i}(\nu)$ and $t_{\mu,i}(\nu)$ 
to denote the starting and terminate times instead of $s_{\mu,i}$ 
	   and $t_{\mu,i}$, 
since they strongly depend on constraint conditions influenced by fiducial machine $\nu$. 
\item[Step 2:] {(Forward Scheduling) With respect to machine 
	   $\mu(\nu<\mu\le M)$  of} the {processing system}, starting and 
	   terminate times of this machine for job $i$ are uniquely 
	   determined by the following relations:
\bea
s_{\mu,i}(\nu)\eq\max\left(
t_{\mu,i-1}(\nu),
t_{\mu-1,i}(\nu)
\right),\\
t_{\mu,i}(\nu)\eq s_{\mu,i}(\nu)+x_{\mu,i}.
\eea
In Step 2, 
$
s_{\nu+1,1}(\nu),t_{\nu+1,1}(\nu),\cdots,
s_{\nu+1,N}(\nu),t_{\nu+1,N}(\nu),
$ $\cdots,
s_{M,1}(\nu),t_{M,1}(\nu),$ $\cdots,
s_{M,N}(\nu),t_{M,N}(\nu)$ are uniquely derived in sequence.

\item[Step 3:]{(Backward Scheduling) With respect to machine 
	   $\mu(1\le\mu<\nu)$  of} the 
	   {processing system}, the starting and 
	   terminate times of this machine for job $i$ are uniquely 
	   determined by the following relations:
\bea
t_{\mu,i}(\nu)\eq\min\left(s_{\mu,i+1}(\nu),s_{\mu+1,i}(\nu)\right),\\
s_{\mu,i}(\nu)\eq t_{\mu,i}(\nu)-x_{\mu,i}.
\eea
In Step 3, 
$
s_{1,1}(\nu),t_{1,1}(\nu),\cdots,
s_{1,N}(\nu),t_{1,N}(\nu),\cdots,
$ $
s_{\nu-1,1}(\nu),t_{\nu-1,1}(\nu),\cdots,$ $
s_{\nu-1,N}(\nu),t_{\nu-1,N}(\nu)
$ are easily determined.

\item[Step 4:]We estimate the makespan of fiducial machine $\nu$ as follows:
\bea
\makespan(\nu)\eq t_{M,N}(\nu)-s_{1,1}(\nu).
\eea
\item[Step 5:]If $\nu<M$, we replace $\nu$ by $\nu+1$ and return to Step 1, otherwise 
this algorithm stops.
\end{description}
As mentioned above, when the processing time table 
$X=\left\{x_{\mu,i}\right\}\in{\bf R}^{M\times N}$ is given 
as a function of the position of the fiducial machine in the {processing system}, the makespan is uniquely 
determined, since the processing time is independently and identically 
distributed according to the probability distribution ${\rm Pr}(x)$.
Thus, we need to average the makespan over multiple randomly generated processing time tables $X$, that is, estimate ${\rm E}_X[\makespan(\nu)]$, in 
order to examine the typical behavior of a particular statistic
measure of  the processing performance of the {processing system.}

\begin{figure}[tb]
\begin{center}
\includegraphics[width=0.9\hsize]{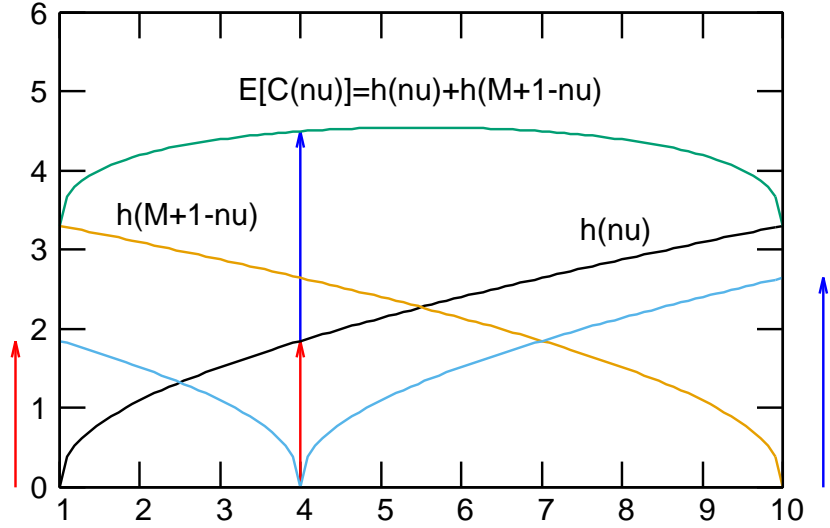}
\caption{\label{shape-function}${\rm E}_X[\makespan(\nu)]=h(\nu)+h(M+1-\nu)$ 
 where $\l={\rm E}_X[X_{\mu,i}]=0$ and $h(\nu)\propto\sqrt{\nu-1}$.
{The horizontal axis shows the position of the fiducial machine and 
the vertical axis shows the makespan or shape function.}}
\end{center}
\end{figure}

Two points should be noticed here. First, from the definitions of forward scheduling and 
backward scheduling, because of system symmetry, 
${\rm E}_X[\makespan(1)]={\rm E}_X[\makespan(M)]$, ${\rm E}_X[\makespan(2)]={\rm E}_X[\makespan(M-1)]$, and ${\rm E}_X[\makespan(3)]={\rm E}_X[\makespan(M-2)]$ are symmetrically and statistically  satisfied, that is, 
\bea
{\rm E}_X[\makespan(\nu)]\eq {\rm E}_X[\makespan(M+1-\nu)].
\eea
{Next} we define the expected terminate time function $h(\nu)$ as
\bea
h(\nu)\eq {\rm E}_X[t_{\nu,N}(\nu)-s_{1,1}(\nu)].
\eea
Then 
\bea
\makespan(\nu)\eq (t_{M,N}(\nu)-s_{\nu,1}(\nu))
-(t_{\nu,N}(\nu)-s_{\nu,1}(\nu))\nn
&&
+(t_{\nu,N}(\nu)-s_{1,1}(\nu)),
\eea
is obtained, and the expectation of makespan can be decomposed into the sum of two terms as follows:
\bea
{\rm E}_X[\makespan(\nu)]\eq h(\nu)+h(M+1-\nu)
-N\l,\qquad\label{eq17}
\eea
where $\l={\rm E}_X[X_{\mu,i}]$. 
Here we replace $L-t_{M,N}(\nu)=s'_{1,1}(M+1-\nu)$ and 
$L-s_{\nu,1}(\nu)=t'_{M+1-\nu,N}(M+1-\nu)$ using an appropriate constant $L$, and
 ${\rm E}_X[t_{M,N}(\nu)-s_{\nu,1}(\nu)]
=
{\rm E}_X[
t'_{M+1-\nu,N}(M+1-\nu)-
s'_{1,1}(M+1-\nu)
]=h(M+1-\nu)
$ is obtained by symmetry. The decomposition of ${\rm E}_X[\makespan(\mu)]$ in \Sref{eq17} is sketched in Fig. 
\ref{shape-function}.

If $M$ and $N$ are sufficiently large, we can prove some properties of 
$h(\nu)$, such as concavity, regardless of the probability distribution of 
$x_{\mu,i}$. Furthermore, 
an explicit form of function $h(\nu)$ can be obtained for some classes of 
probability distributions. We will present these results in detail in section \ref{mathematical-discussion}.

\section{Numerical Experiments\label{sec3}}

\subsection{Exponential Distribution}
\begin{figure}[t]
\begin{center}
\includegraphics[width=0.8\hsize]{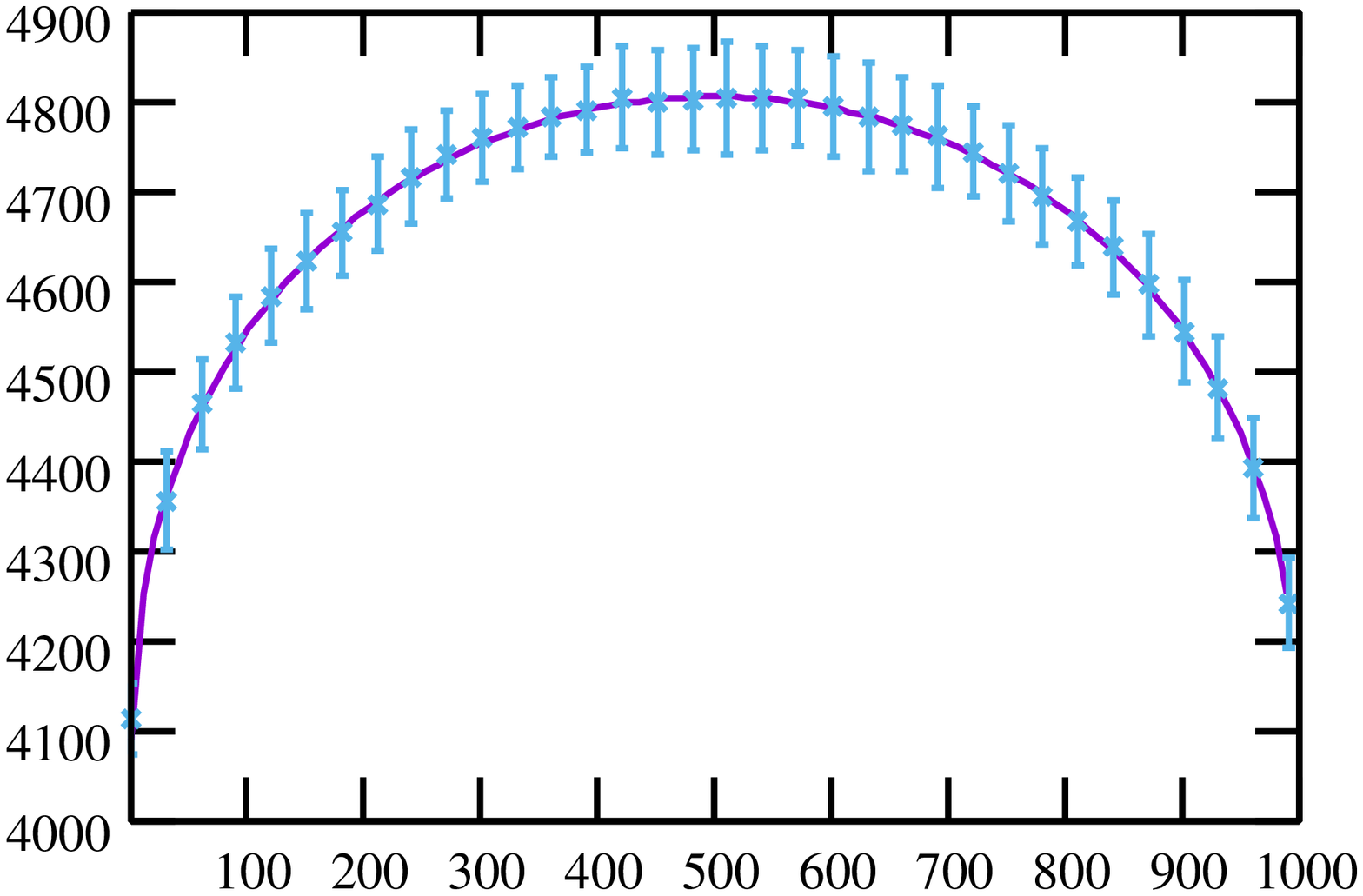} 
\caption{\label{fig-exp-200}$M=1000$ and $N=200$, exponential distribution.
{The horizontal axis shows the position of the fiducial machine 
 $\mu$ and the vertical axis shows the makespan $\makespan(\mu)$.}
}
\vspace{0.5cm}
\includegraphics[width=0.8\hsize]{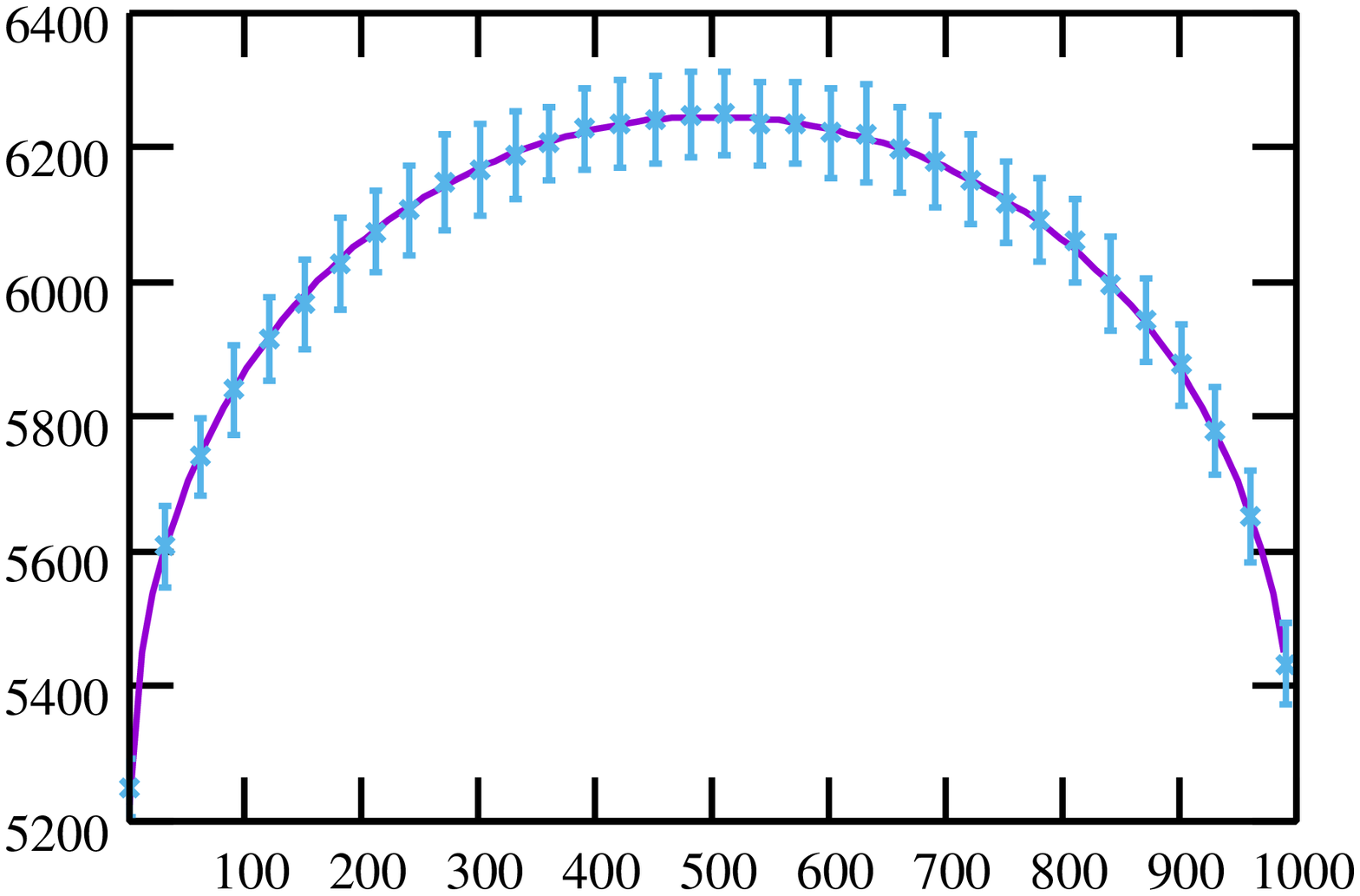} 
\caption{\label{fig-exp-400}$M=1000$ and $N=400$, exponential distribution.
{The horizontal axis shows the position of the fiducial machine 
 $\mu$ and the vertical axis shows the makespan $\makespan(\mu)$.}
}
\end{center}
\end{figure}
\begin{figure}[t]
\begin{center}
\includegraphics[width=0.8\hsize]{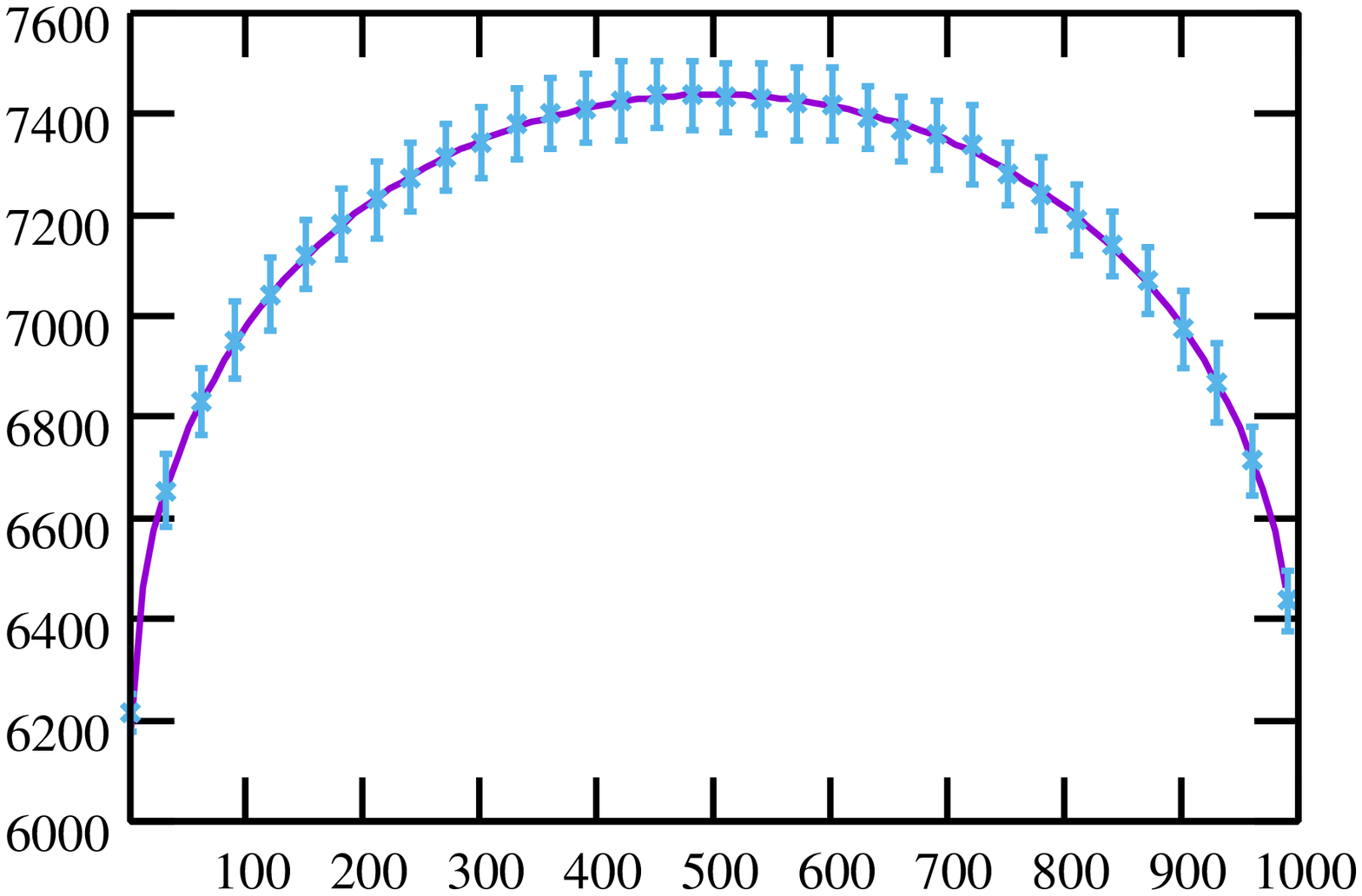} 
\caption{\label{fig-exp-600}$M=1000$ and $N=600$, exponential distribution.
{The horizontal axis shows the position of the fiducial machine 
 $\mu$ and the vertical axis shows the makespan $\makespan(\mu)$.}
}
\vspace{0.5cm}
\includegraphics[width=0.8\hsize]{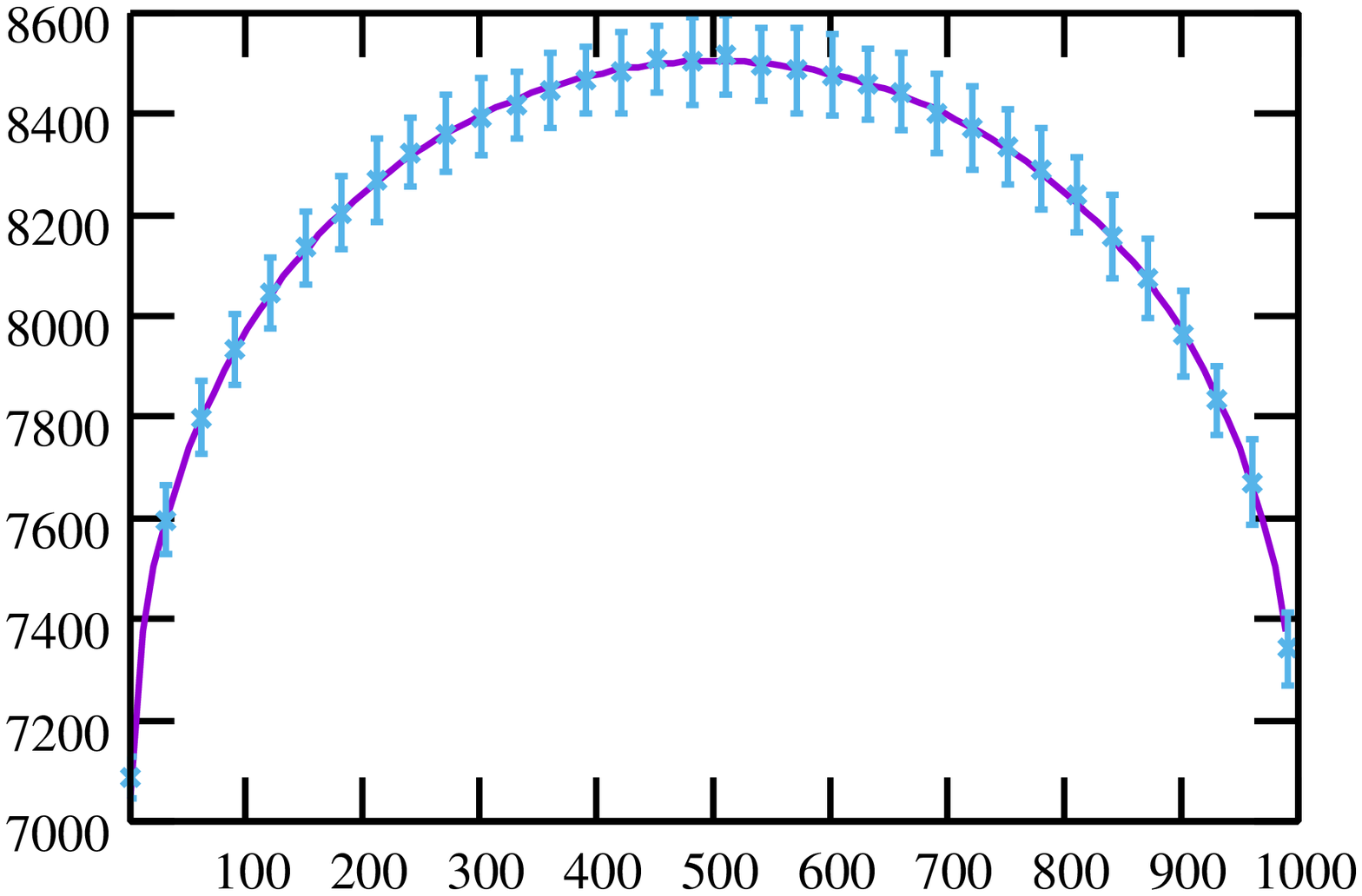} 
\caption{\label{fig-exp-800}$M=1000$ and $N=800$, exponential distribution.
{The horizontal axis shows the position of the fiducial machine 
 $\mu$ and the vertical axis shows the makespan $\makespan(\mu)$.}
}
\end{center}
\end{figure}

\begin{table}[t]
\begin{center}
\begin{tabular}{|c|c|c|c|c|l|}
\hline
$M$&$N$&$A$&$B$&$\a$&Figure\\
\hline
$1000$&$200$
&$53.1185$
&$2333.8463$
&$0.50652$
&Fig.\ref{fig-exp-200}
\\
\hline
$1000$&$400$
&$72.8693$
&$2660.1637$
&$0.51548$
&Fig.\ref{fig-exp-400}
\\
\hline
$1000$&$600$
&$96.5423$
&$3107.2000$
&$0.50059$
&Fig.\ref{fig-exp-600}
\\
\hline
$1000$&$800$
&$106.8016$
&$3443.8947$
&$0.50941$
&Fig.\ref{fig-exp-800}
\\
\hline
\end{tabular}
\caption{Parameter estimates for fitting functions (plotted in the figures) for the case that the processing time is exponentially distributed. \label{tab-exp}}
\end{center}
\end{table}
\clearpage

In Fig. \ref{shape-function}, an example of the expected terminate time $h(\nu)$ is depicted 
as a concave function. 
We can prove that indeed if $M$ and $N$ are sufficiently large, 
$h(\nu)$ will be a concave function regardless of the probability 
distribution of $X_{\mu,i}$. 
However, we will first discuss the relationship between the expectation of makespan 
and the position of the fiducial machine in the {processing system} when the processing time is independently and identically distributed according to an exponential distribution.
In this discussion, the following representation of the exponential distribution is used:
\bea
f(X_{\mu,i})\eq
\left\{
\begin{array}{ll}
\f{1}{\l}e^{-\f{X_{\mu,i}}{\l}}&X_{\mu,i}>0\\
0&\text{otherwise}
\end{array}
\right.,
\eea
where $\l$ is the scale parameter and is such that
${\rm E}_X[X_{\mu,i}]=\l$ and 
${\rm Var}[X_{\mu,i}]=
{\rm E}_X[X^2_{\mu,i}]
-({\rm E}_X[X_{\mu,i}])^2
=
\l^2$.

Here we apply $\l=2$ in numerical experiments, the number of machines in 
{this system} is $M=1000$, and the number of jobs processed in 
the {processing system} is $N=200,400,600$, or $800$. These results for the 
average of makespan with respect to the index of fiducial machine $\nu$ 
are shown in {Figs. \ref{fig-exp-200} 
to \ref{fig-exp-800},} respectively. The horizontal and vertical axes in these figures indicate the index of the fiducial machine $\nu$ and the average of the makespan ${\rm E}_X[\makespan(\nu)]$ as evaluated by numerical experiments, respectively.  
Specifically, the symbols and error bars are results from evaluating $100$ random 
processing time tables.  As shown, 
the behavior of the mean makespan depends not on the position of the 
bottleneck machine but on the position of the fiducial machine in {this 
system} because the bottleneck machine equivalently exists in the 
{processing system}. Moreover using equation (\ref{eq17}), the fitting function (solid lines in these figures) is assumed as 
\bea
{\rm E}_X[\makespan(\nu)]\eq A(\nu-1)^\a+A(M-\nu)^\a+B,\qquad\label{eq21}
\eea
where 
$B=2b+c(M-1)-N{\rm E}_X[X_{\nu,i}]$. This equation corresponds to the \shape functions $h(\nu)=A(\nu-1)^\a+b+c(\nu-1)$ and 
$h(M+1-\nu)=A(M-\nu)^\a+b+c(M-\nu)$ already used. These figures show that the results of numerical experiments and fitting functions are 
consistent with each other. These estimates of parameters are listed in table 
\ref{tab-exp}. From this table, it is expected that  $\a$ is constant 
with respect to the number of jobs in the case of this exponential distribution; that is, the assumption in equation 
(\ref{eq17}) that the mean of the makespan can be expressed by two \shape functions is satisfied.

\begin{figure}[t] 
\begin{center}
\includegraphics[width=0.8\hsize]{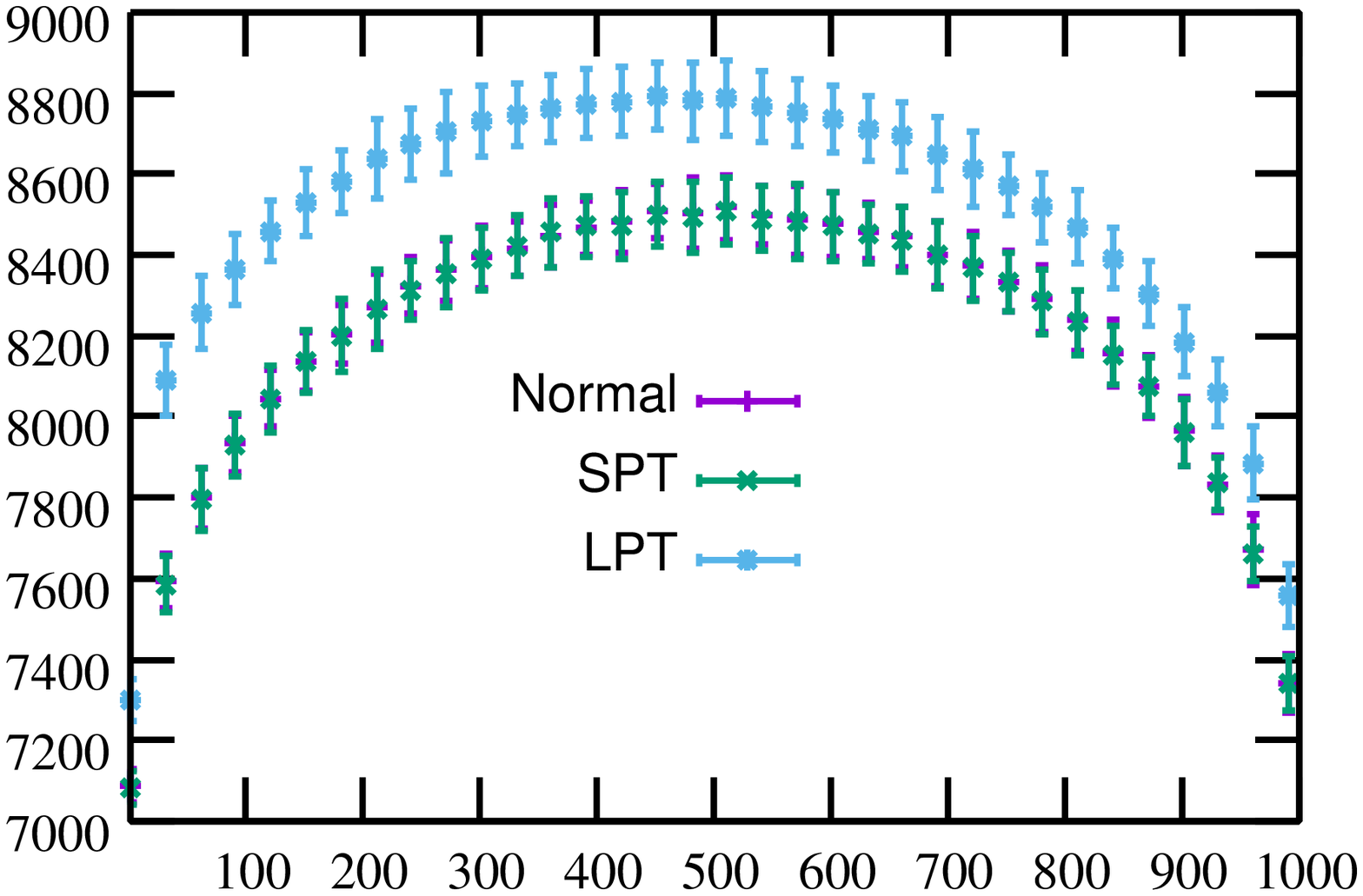} 
\caption{\label{dispatching-rule-ohsawa0}Makespans under normal, SPT, and LPT rules. $M=1000$ and $N=800$,  exponential 
 distribution with mean $\l=2$. 
{The horizontal axis shows the position of the fiducial machine 
 $\mu$ and the vertical axis shows the makespan $\makespan(\mu)$.}
}
\vspace{0.5cm}
\includegraphics[width=0.8\hsize]{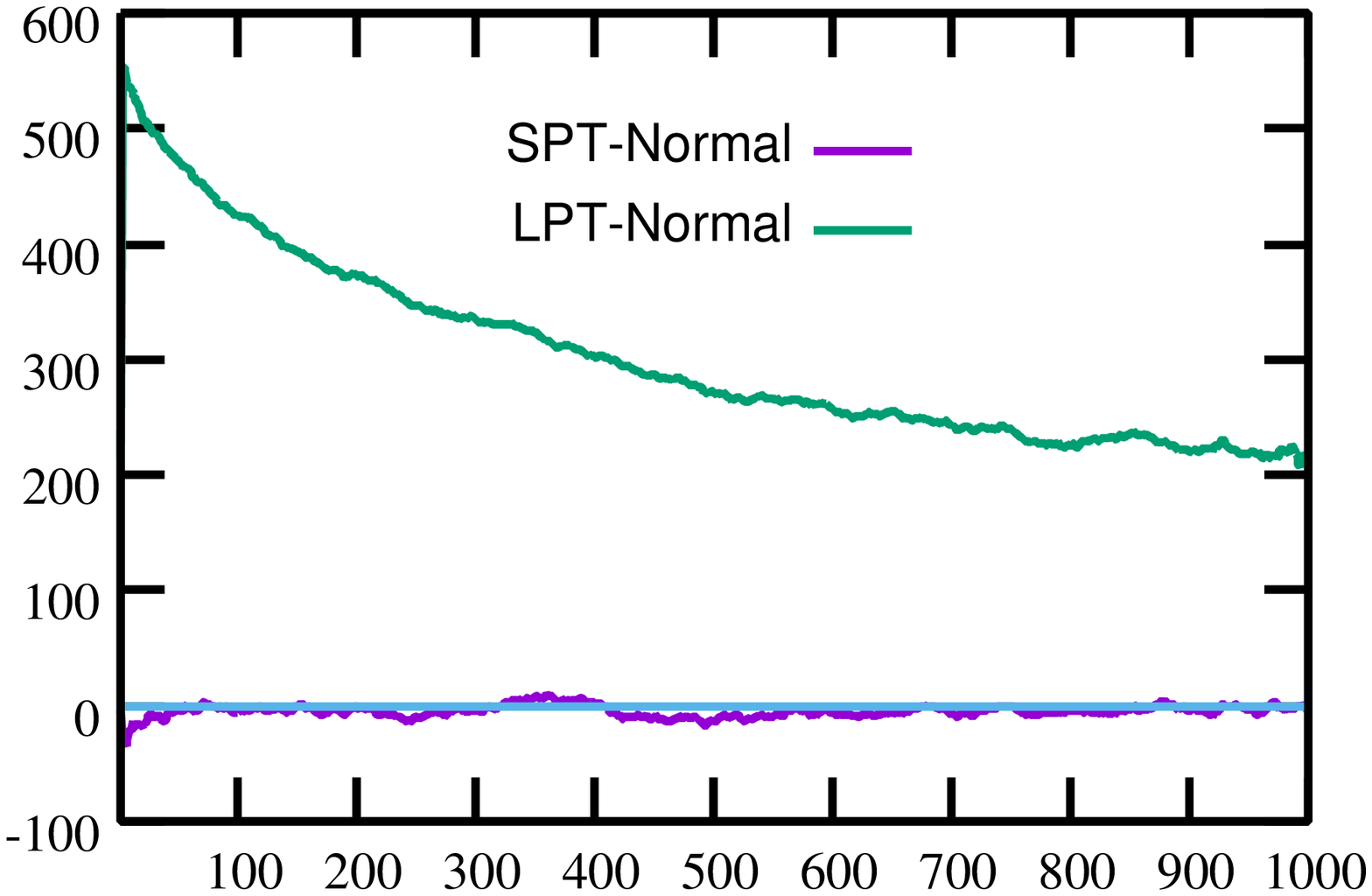} 
\caption{\label{dispatching-rule-del-ohsawa0}Differences in makespan as plotted in {Fig.} \ref{dispatching-rule-ohsawa0} between normal and  dispatching rules.
{The horizontal axis shows the position of the fiducial machine 
 $\mu$ and the vertical axis shows the makespan $\makespan(\mu)$.}
}
\end{center}
\end{figure}

Lastly, we compare the above flowshop scheduling  results with those using two of the most discussed 
and applied dispatching rules, shortest processing time (SPT) and 
longest processing time (LPT). For the case $N=800$ and $M=1000$, the results of applying the two dispatching rules are given in terms of the processing times of machine 
$1$ in {Figs. \ref{dispatching-rule-ohsawa0} and 
 \ref{dispatching-rule-del-ohsawa0}.} As shown, the behaviors of the expectation of makespan for the normal rule 
scheduling  and those of SPT rule scheduling are statistically consistent with each other because the processing time is independently and identically distributed between jobs and machines
in 
the case of large systems, which results in attempts at optimization being 
{canceled out}  by the randomness of processing time, causing the expected effect 
of {these dispatching rules}  to be negligible. In fact, $\mathop{\min}_\nu({\rm E}_X[\makespan_{SPT}(\nu)]-{\rm E}_X[\makespan(\nu)])/{\rm E}_X[\makespan(\nu)]\simeq-0.0047$ and $\mathop{\max}_\nu({\rm E}_X[\makespan_{SPT}(\nu)]-{\rm E}_X[\makespan(\nu)])/{\rm E}_X[\makespan(\nu)]\simeq0.0012$, where $\makespan_{SPT}(\nu)$ and $\makespan(\nu)$ are the makespans under the SPT rule and the normal rule, respectively, so that the expected improvement effect is only $0.47\%$, and thus intuitively it is not cost-effective to try to obtain an improvement using the dispatching rule.  
{In contrast}, as shown in {Figs. \ref{dispatching-rule-ohsawa0} and \ref{dispatching-rule-del-ohsawa0},} the makespan under LPT scheduling 
{which is not optimal} is larger than that under normal scheduling. However, $\mathop{\min}_\nu({\rm E}_X[\makespan_{LPT}(\nu)]-{\rm E}_X[\makespan(\nu)])/{\rm E}_X[\makespan(\nu)]\simeq0.027$ and $\mathop{\max}_\nu({\rm E}_X[\makespan_{LPT}(\nu)]-{\rm E}_X[\makespan(\nu)])/{\rm E}_X[\makespan(\nu)]\simeq0.077$, where $\makespan_{LPT}(\nu)$ is the makespan under the LPT rule, 
{so LPT is at most $7.7\%$  worse than normal scheduling}.

\subsection{Discrete Uniform Distribution}

\begin{figure}[t]
\begin{center}
\includegraphics[width=0.8\hsize]{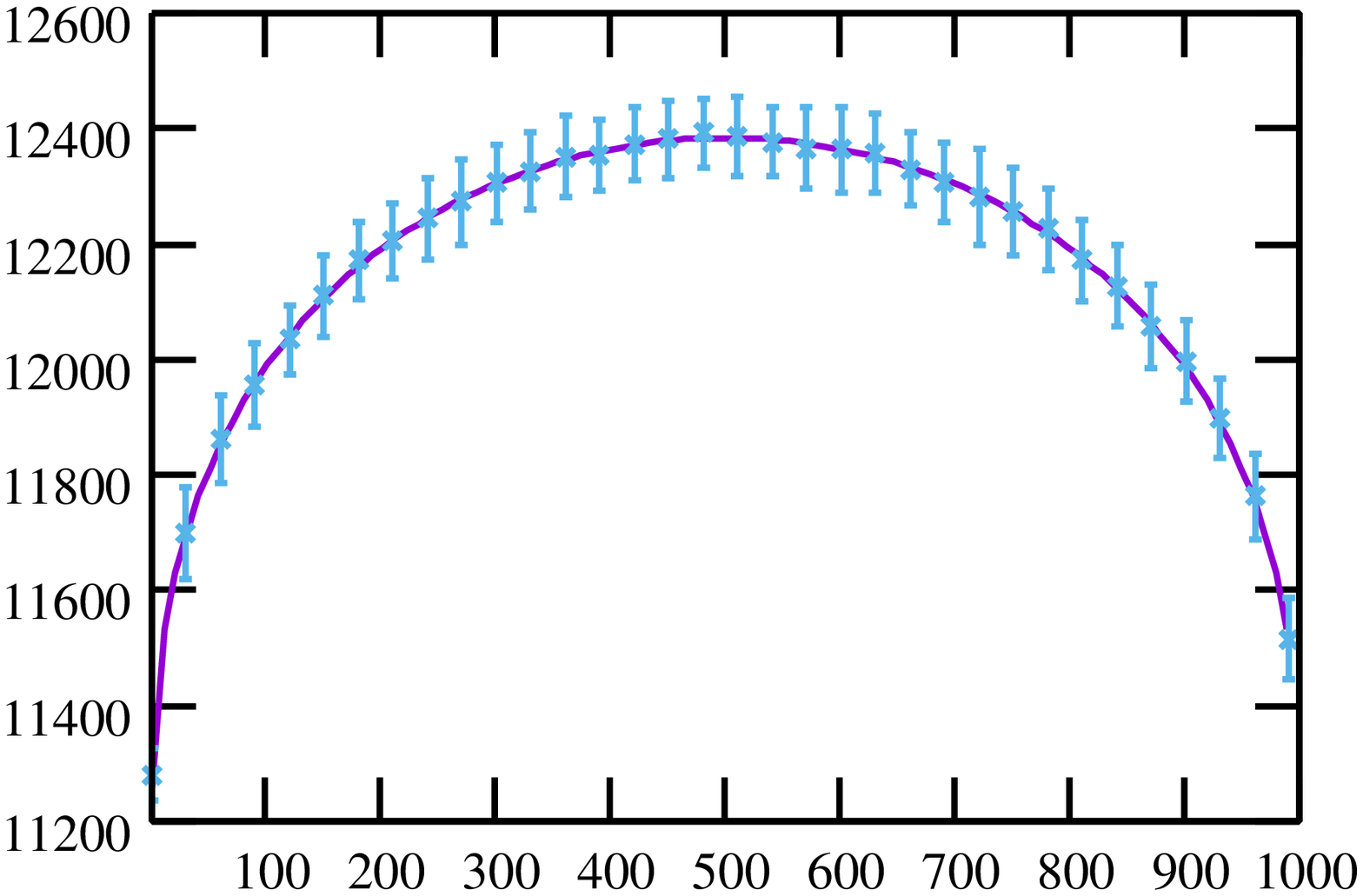} 
\caption{\label{fig-dis-uni-200}$M=1000$ and $N=200$,  discrete uniform distribution.
{The horizontal axis shows the position of the fiducial machine 
 $\mu$ and the vertical axis shows the makespan $\makespan(\mu)$.}
}
\vspace{0.5cm}
\includegraphics[width=0.8\hsize]{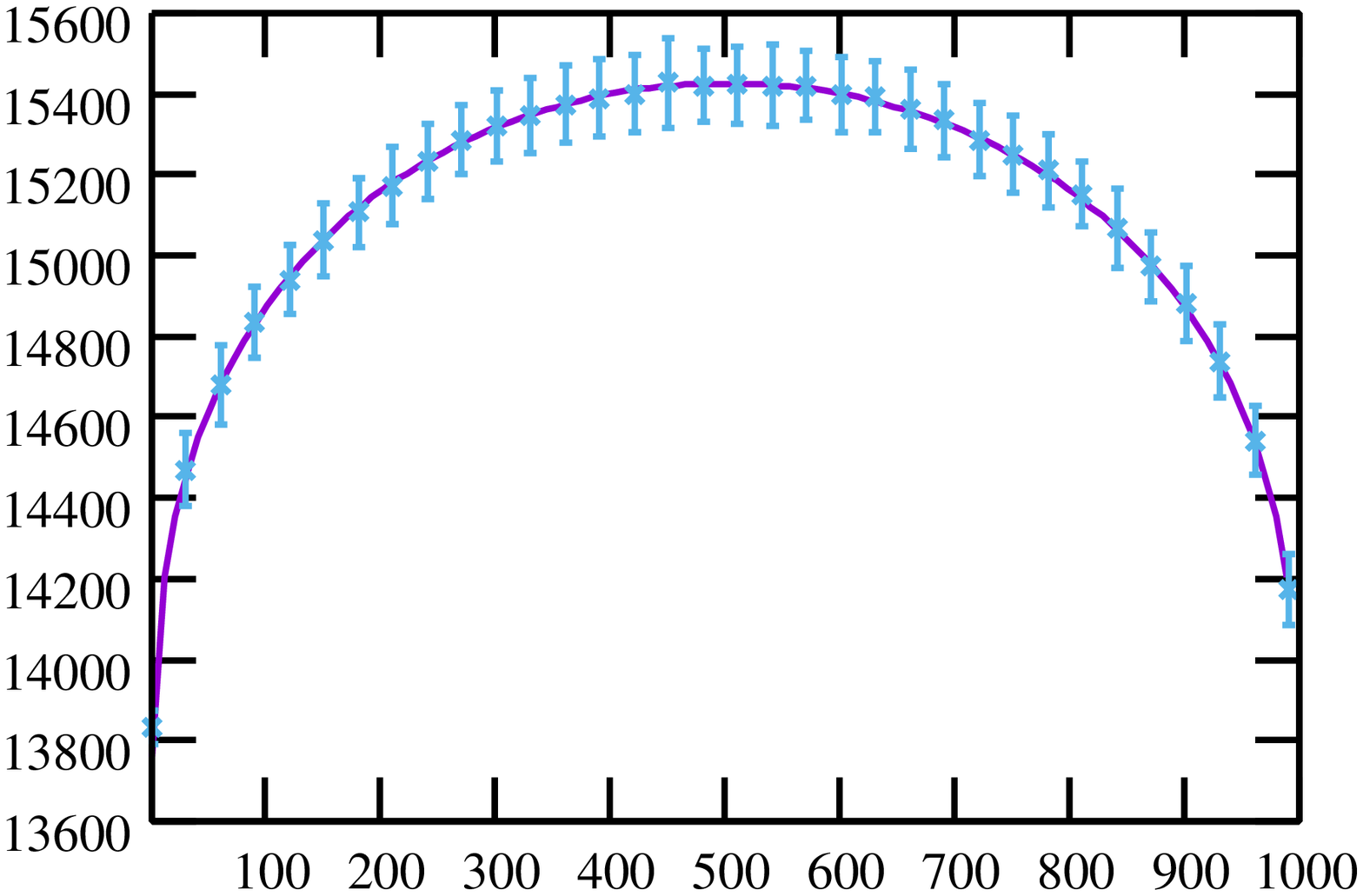} 
\caption{\label{fig-dis-uni-400}$M=1000$ and $N=400$,  discrete uniform distribution.
{The horizontal axis shows the position of the fiducial machine 
 $\mu$ and the vertical axis shows the makespan $\makespan(\mu)$.}
}
\end{center}
\end{figure}
\begin{figure}[t]
\begin{center}
\includegraphics[width=0.8\hsize]{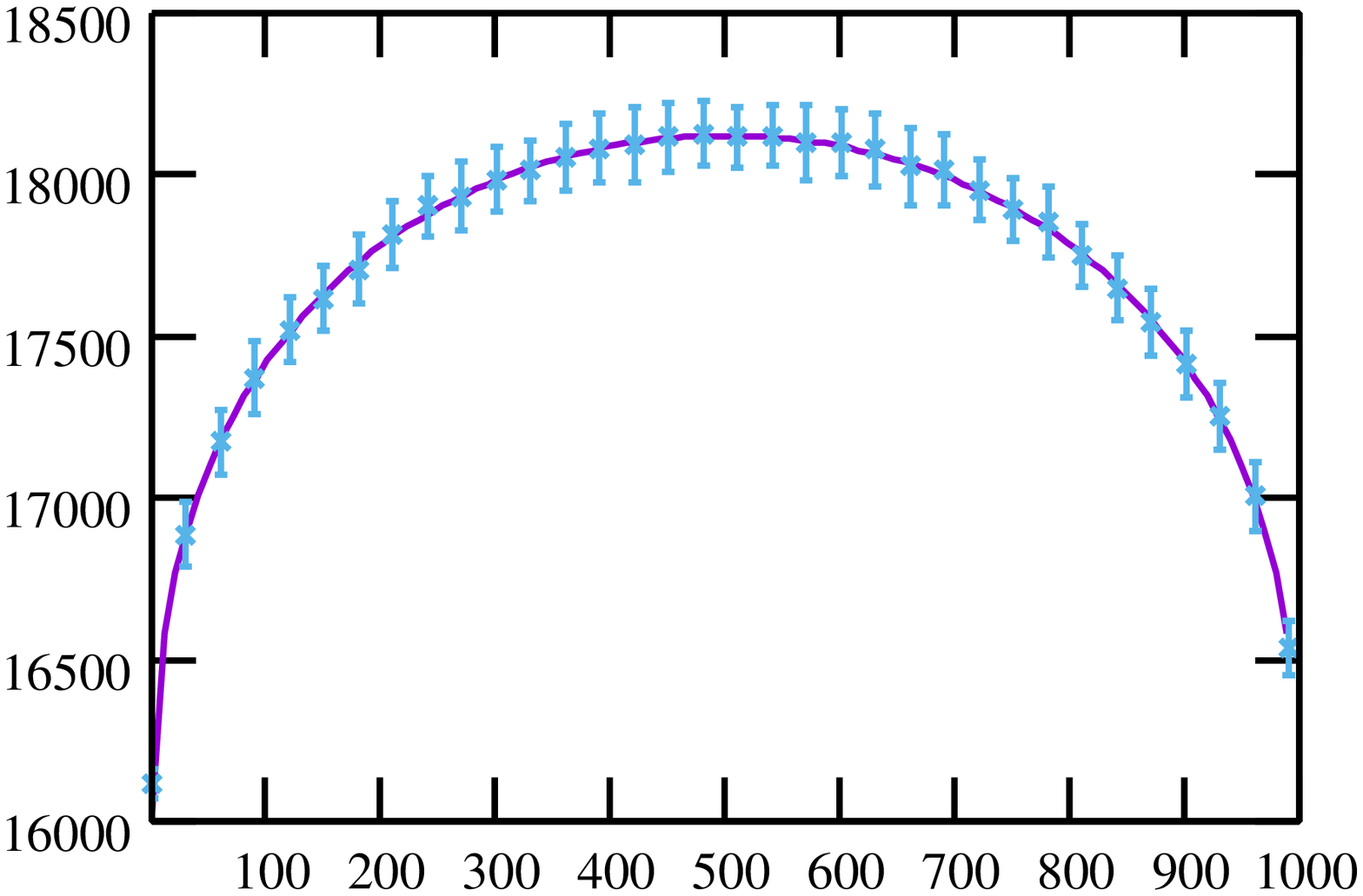} 
\caption{\label{fig-dis-uni-600}$M=1000$ and $N=600$,  discrete uniform distribution.
{The horizontal axis shows the position of the fiducial machine 
 $\mu$ and the vertical axis shows the makespan $\makespan(\mu)$.}
}
\vspace{0.5cm}
\includegraphics[width=0.8\hsize]{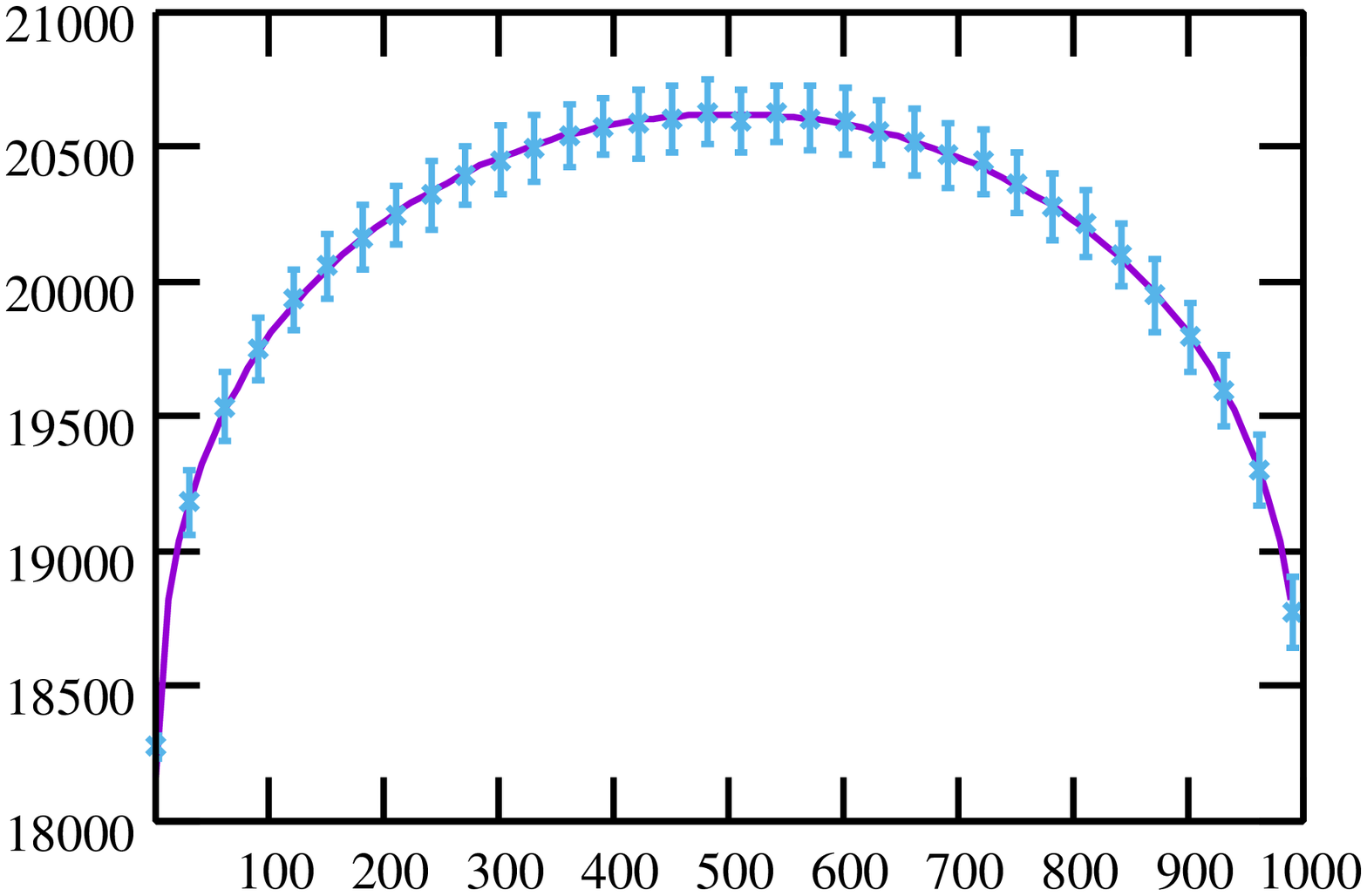} 
\caption{\label{fig-dis-uni-800}$M=1000$ and $N=800$,  discrete uniform distribution.
{The horizontal axis shows the position of the fiducial machine 
 $\mu$ and the vertical axis shows the makespan $\makespan(\mu)$.}
}
\end{center}
\end{figure}

\begin{table}[t]
\begin{center}
\begin{tabular}{|c|c|c|c|c|l|}
\hline
$M$&$N$&$A$&$B$&$\a$&Figure\\
\hline
$1000$&$200$
&$93.7442$
&$8661.74$
&$0.48092$
&Fig.\ref{fig-dis-uni-200}
\\
\hline
$1000$&$400$
&$164.7893$
&$10218.35$
&$0.44421$
&Fig.\ref{fig-dis-uni-400}
\\
\hline
$1000$&$600$
&$202.0881$
&$11521.40$
&$0.44950$
&Fig.\ref{fig-dis-uni-600}
\\
\hline
$1000$&$800$
&$241.3292$
&$12917.00$
&$0.44587$
&Fig.\ref{fig-dis-uni-800}
\\
\hline
\end{tabular}
\caption{Parameter estimates for fitting functions (plotted in the figures) in the case that processing time is distributed according to a discrete uniform distribution.\label{tab-dis-uni}
}
\end{center}
\end{table}
\clearpage

Next, we also discuss the relationship between the expectation of 
makespan and the position of the fiducial machine in {this system} when 
the processing time 
is independently and identically distributed according to a discrete uniform distribution, specifically, the following distribution:
\bea
{\rm Pr}(X_{\mu,i})\eq
\left\{
\begin{array}{ll}
\f{1}{13}&X_{\mu,i}=1,2,\cdots,13\\
0&\text{otherwise}
\end{array}
\right.,\qquad
\eea
which is such that $
{\rm E}_X[X_{\mu,i}]=7$ and 
$
{\rm Var}[X_{\mu,i}]={\rm E}_X[X^2_{\mu,i}]-({\rm E}_X[X_{\mu,i}])^2
=14$. 
As in the previous subsection, the 
number of machines is $M=1000$ and the number of jobs is $N=200,400,600$, 
or $800$. The results are shown in 
Fig. \ref{fig-dis-uni-200} to 
Fig. \ref{fig-dis-uni-800}. As shown, the expectation of makespan under this distribution 
behaves similarly to the case of an exponential distribution. The estimated parameters in 
equation (\ref{eq21}) for the numerical experiments are listed in table \ref{tab-dis-uni}.

\subsection{Continuous Uniform Distribution}
\begin{figure}[t]
\begin{center}
\includegraphics[width=0.8\hsize]{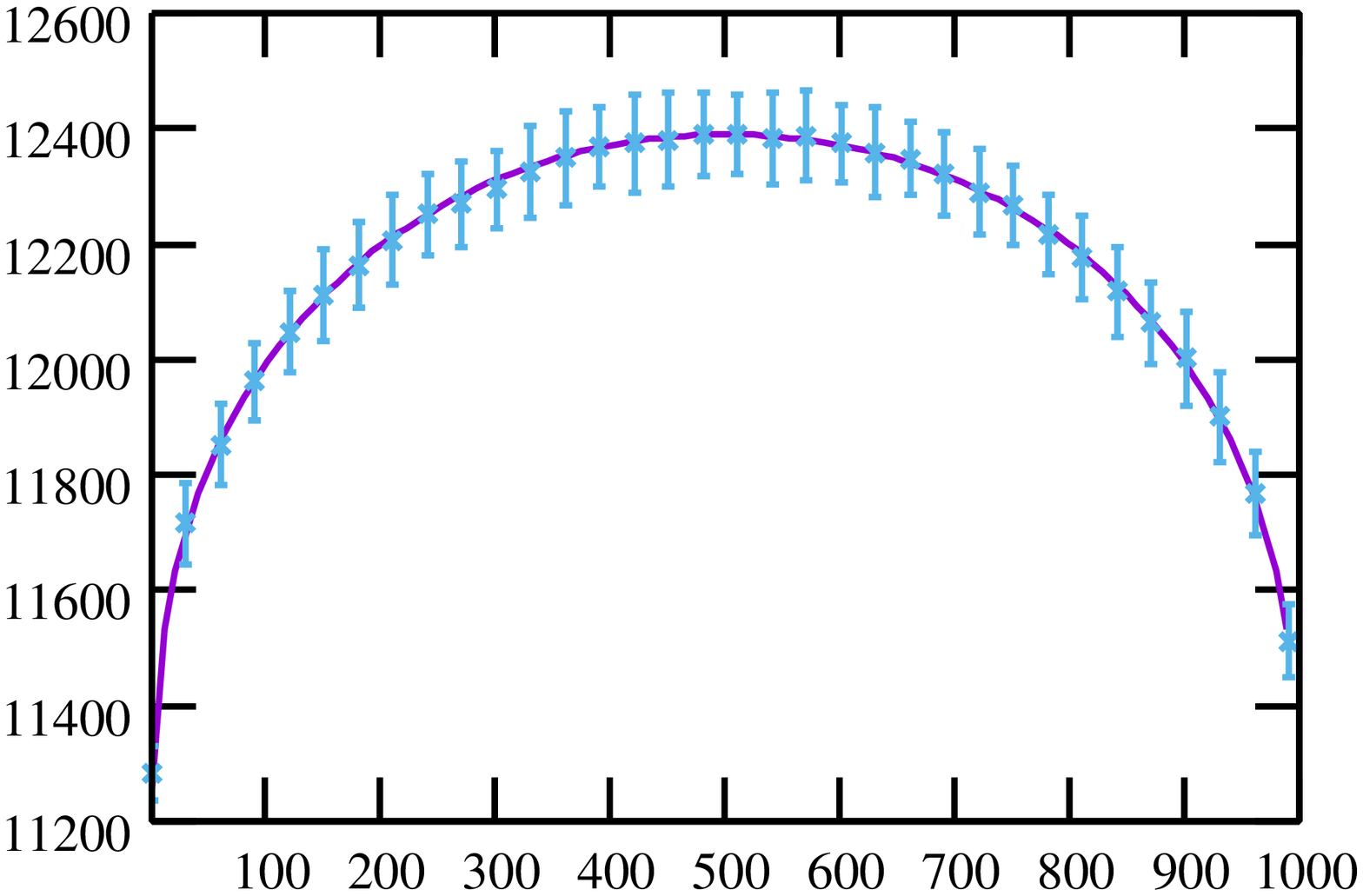} 
\caption{\label{fig-cont-uni-200}$M=1000$ and $N=200$, continuous uniform distribution.{The horizontal axis shows the position of the fiducial machine 
 $\mu$ and the vertical axis shows the makespan $\makespan(\mu)$.}
}
\vspace{0.5cm}
\includegraphics[width=0.8\hsize]{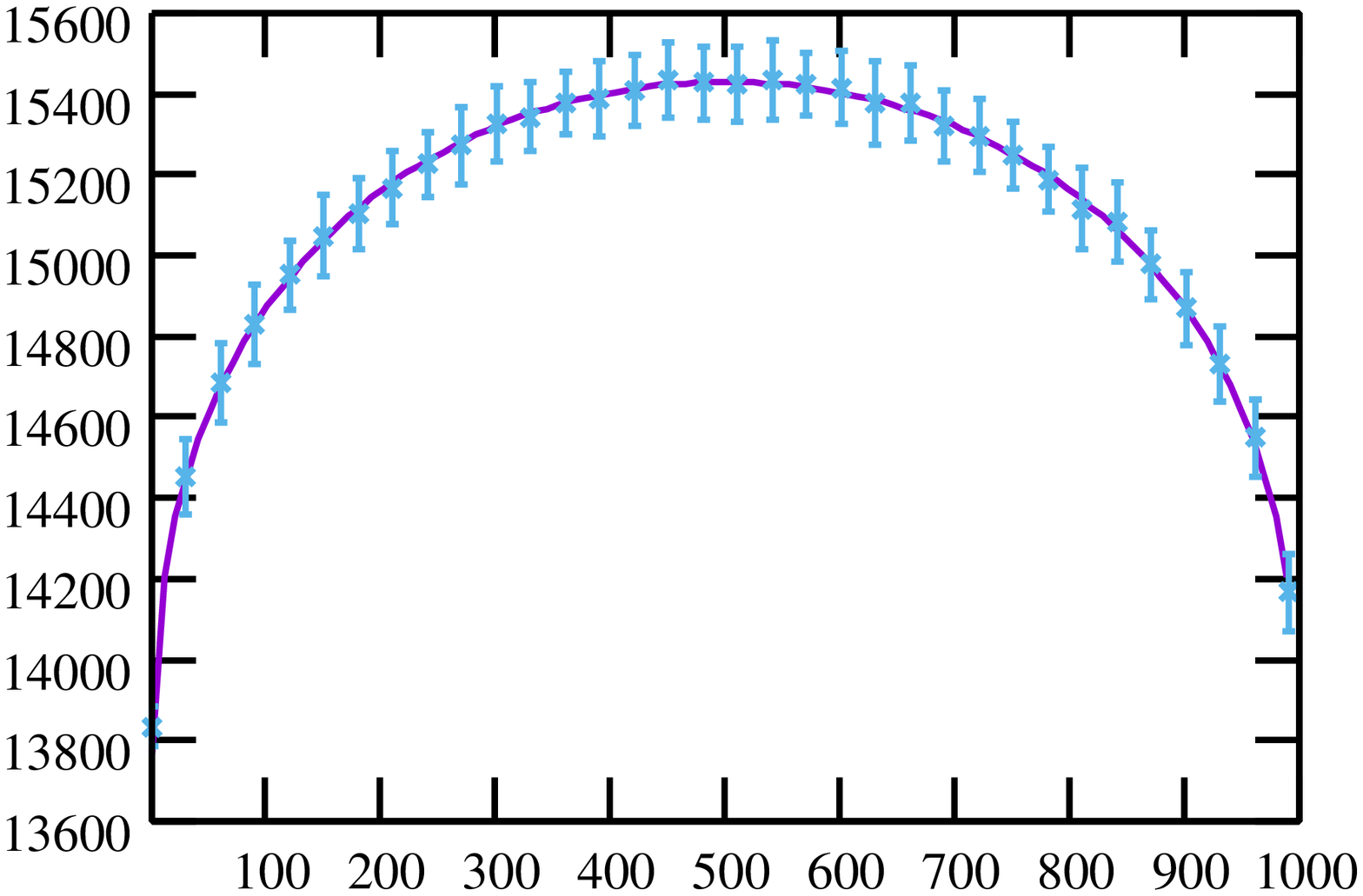} 
\caption{\label{fig-cont-uni-400}$M=1000$ and $N=400$, continuous uniform distribution.{The horizontal axis shows the position of the fiducial machine 
 $\mu$ and the vertical axis shows the makespan $\makespan(\mu)$.}
}
\end{center}
\end{figure}
\begin{figure}[t]
\begin{center}
\includegraphics[width=0.8\hsize]{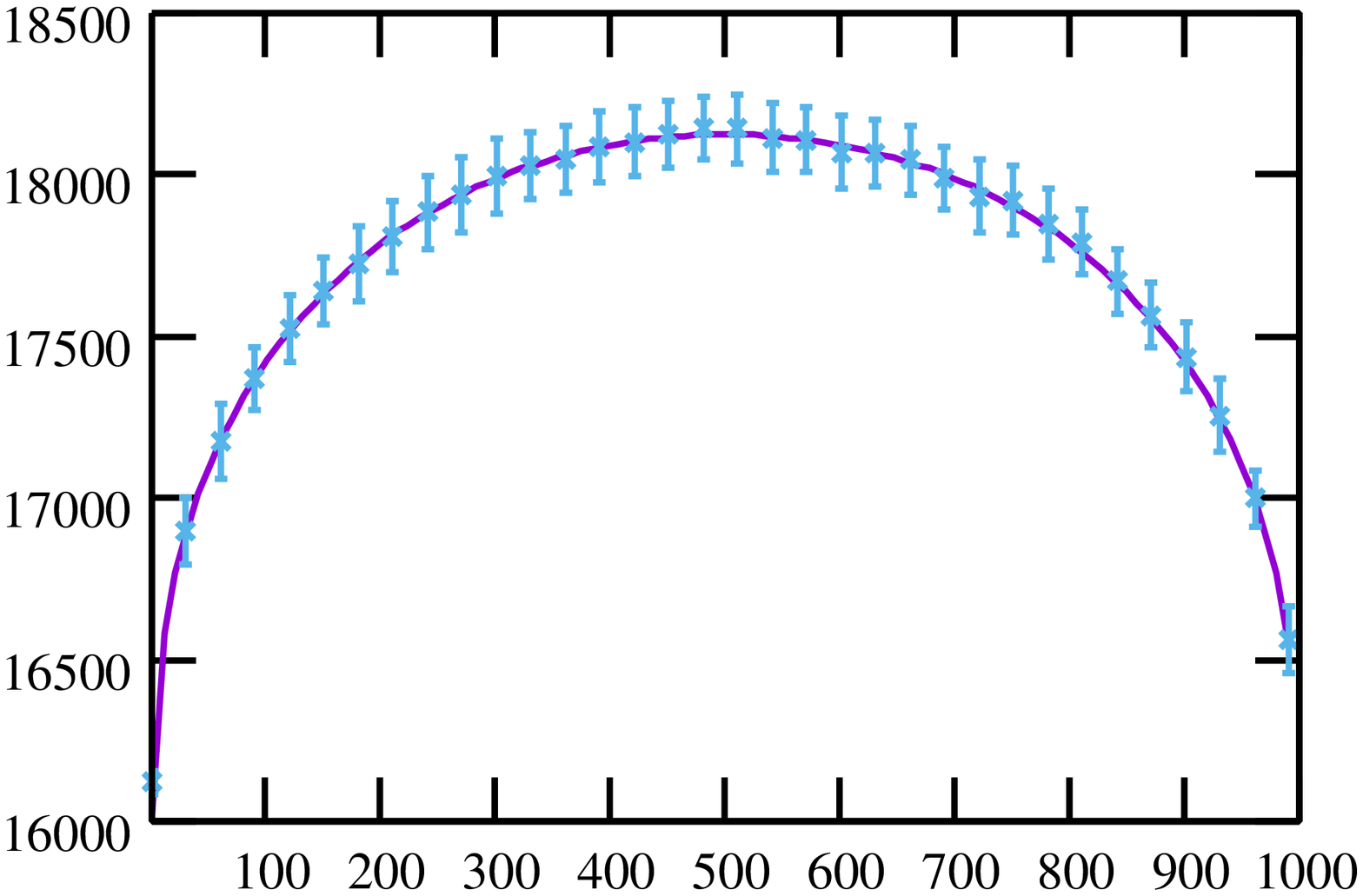} 
\caption{\label{fig-cont-uni-600}$M=1000$ and $N=600$, continuous uniform distribution.{The horizontal axis shows the position of the fiducial machine 
 $\mu$ and the vertical axis shows the makespan $\makespan(\mu)$.}
}
\vspace{0.5cm}
\includegraphics[width=0.8\hsize]{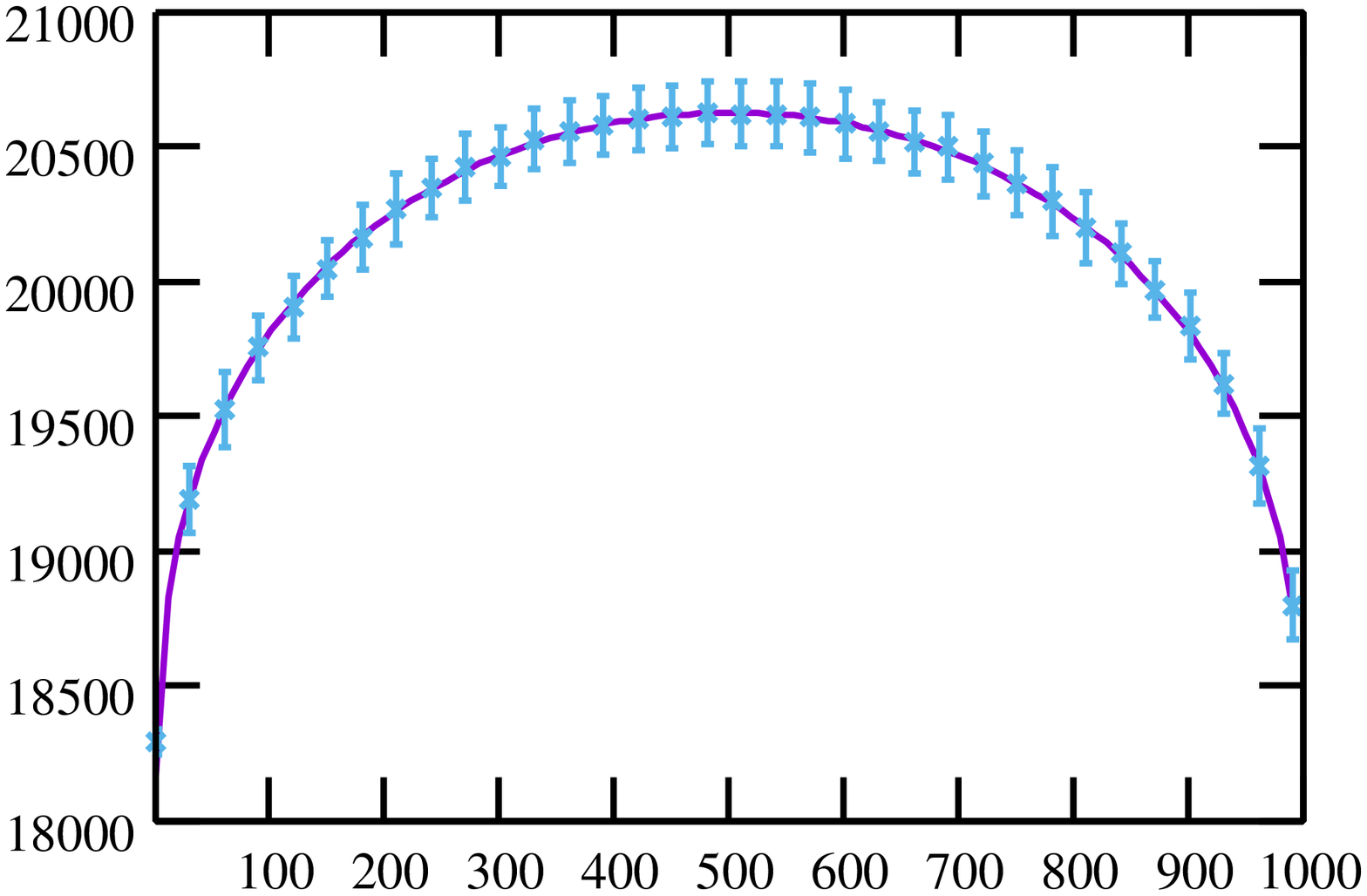} 
\caption{\label{fig-cont-uni-800}$M=1000$ and $N=800$, continuous uniform distribution.{The horizontal axis shows the position of the fiducial machine 
 $\mu$ and the vertical axis shows the makespan $\makespan(\mu)$.}
}

\end{center}
\end{figure}

\begin{table}[t]
\begin{center}
\begin{tabular}{|c|c|c|c|c|l|}
\hline
$M$&$N$&$A$&$B$&$\a$&Figure\\
\hline
$1000$&$200$
&$95.4363$
&$8661.48$
&$0.47832$
&Fig.\ref{fig-cont-uni-200}
\\
\hline
$1000$&$400$
&$159.8063$
&$10180.61$
&$0.45043$
&Fig.\ref{fig-cont-uni-400}
\\
\hline
$1000$&$600$
&$206.0768$
&$11521.64$
&$0.44645$
&Fig.\ref{fig-cont-uni-600}
\\
\hline
$1000$&$800$
&$246.9541$
&$12975.25$
&$0.44105$
&Fig.\ref{fig-cont-uni-800}
\\
\hline
\end{tabular}
\caption{Parameter estimates  for fitting functions (plotted in the figures) for the case that processing time is distributed according to a continuous uniform distribution.\label{tab-cont-uni}}
\end{center}
\end{table}
\clearpage

Next, in order to be able to make a comparison with the previous discrete uniform case, we consider a continuous uniform distribution having the same mean and variance, that is, the density function 
\bea
f(X_{\mu,i})\eq
\left\{
\begin{array}{ll}
\f{1}{2\sqrt{42}}&
\left|X_{\mu,i}-7\right|
\le \sqrt{42}\\
0&\text{otherwise}
\end{array}
\right.,
\eea
is used and is such that ${\rm E}_X[X_{\mu,i}]=7$ and ${\rm Var}[X_{\mu,i}]={\rm E}_X[X^2_{\mu,i}]-({\rm E}_X[X_{\mu,i}])^2=14$. 

Fig. \ref{fig-cont-uni-200} to 
Fig. \ref{fig-cont-uni-800} show plots of the expectations of makespan for the cases of $N=200, 400, 600$, and $800$ with $M=1000$. The parameter estimates are listed in table \ref{tab-cont-uni}. Comparing 
table \ref{tab-dis-uni} and table \ref{tab-cont-uni}, it is clear that 
the behaviors of makespan for the two cases are similar up to a statistical fluctuation.

\subsection{Chi-Squared Distribution}
\begin{figure}[t]
\begin{center}
\includegraphics[width=0.8\hsize]{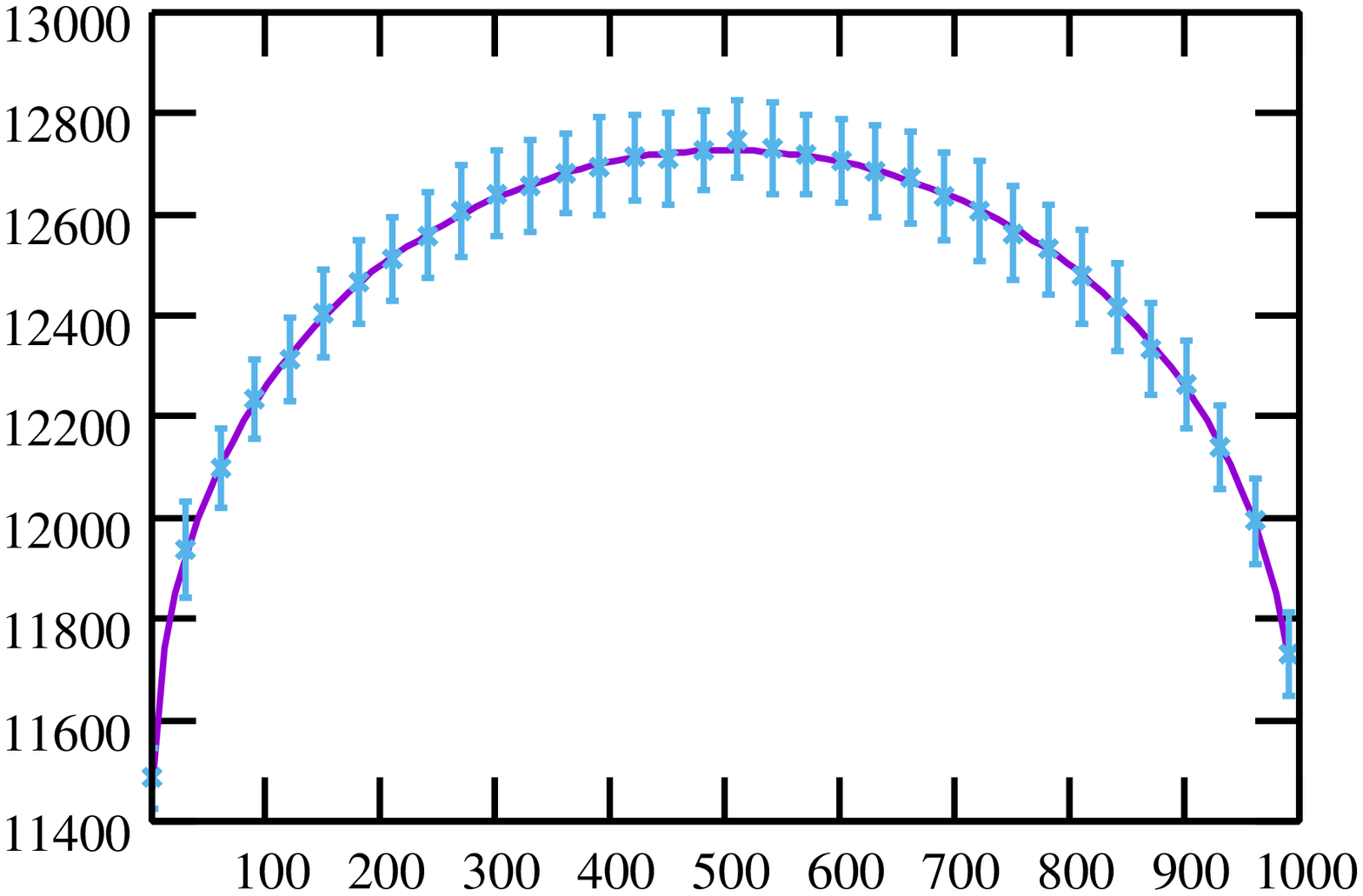} 
\caption{\label{fig-chi-sq-200}$M=1000$ and $N=200$, chi-squared distribution.{The horizontal axis shows the position of the fiducial machine 
 $\mu$ and the vertical axis shows the makespan $\makespan(\mu)$.}
}
\vspace{0.5cm}
\includegraphics[width=0.8\hsize]{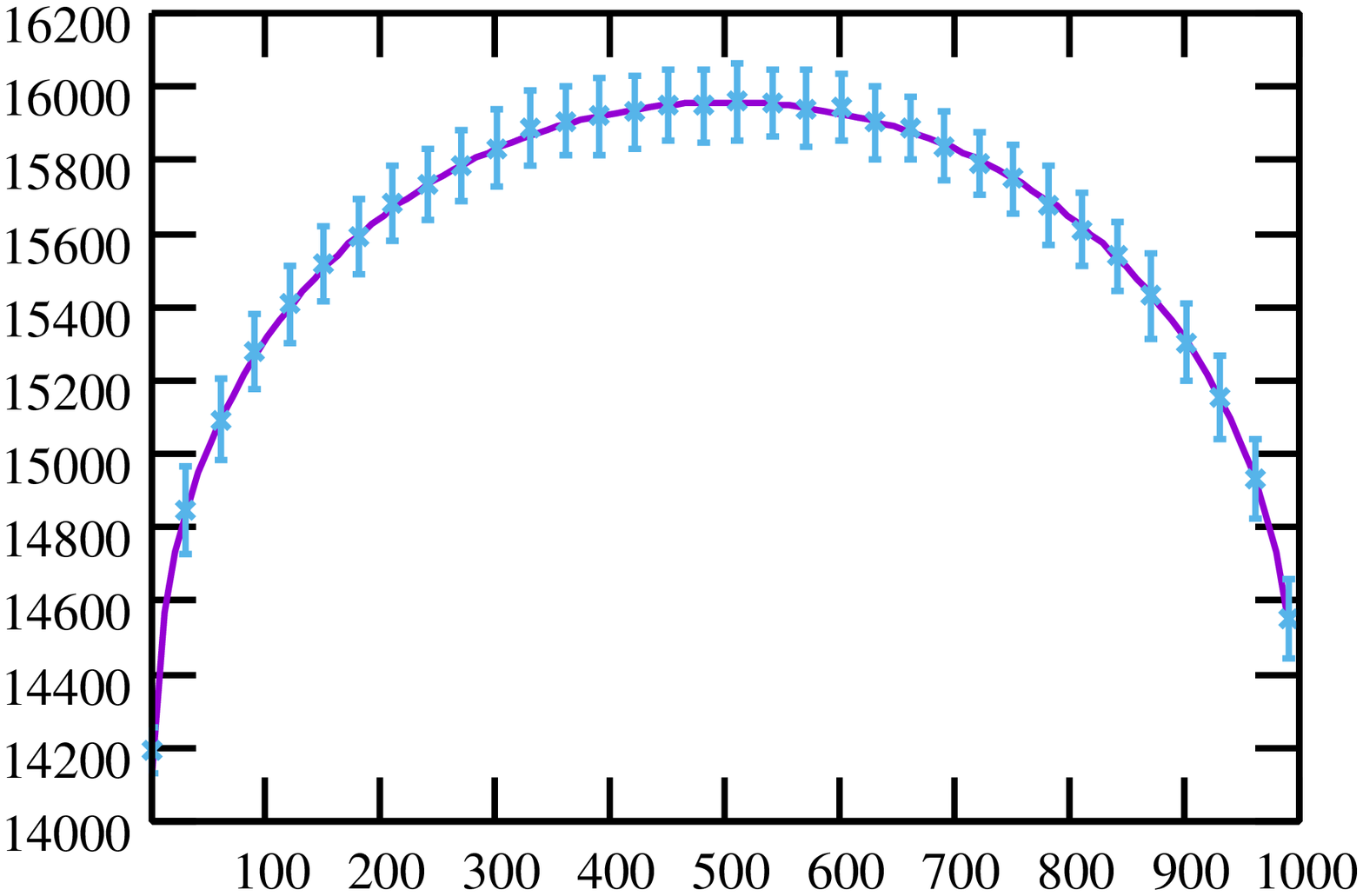} 
\caption{\label{fig-chi-sq-400}$M=1000$ and $N=400$, chi-squared distribution.{The horizontal axis shows the position of the fiducial machine 
 $\mu$ and the vertical axis shows the makespan $\makespan(\mu)$.}
}
\end{center}
\end{figure}

\begin{figure}[t]
\begin{center}
\includegraphics[width=0.8\hsize]{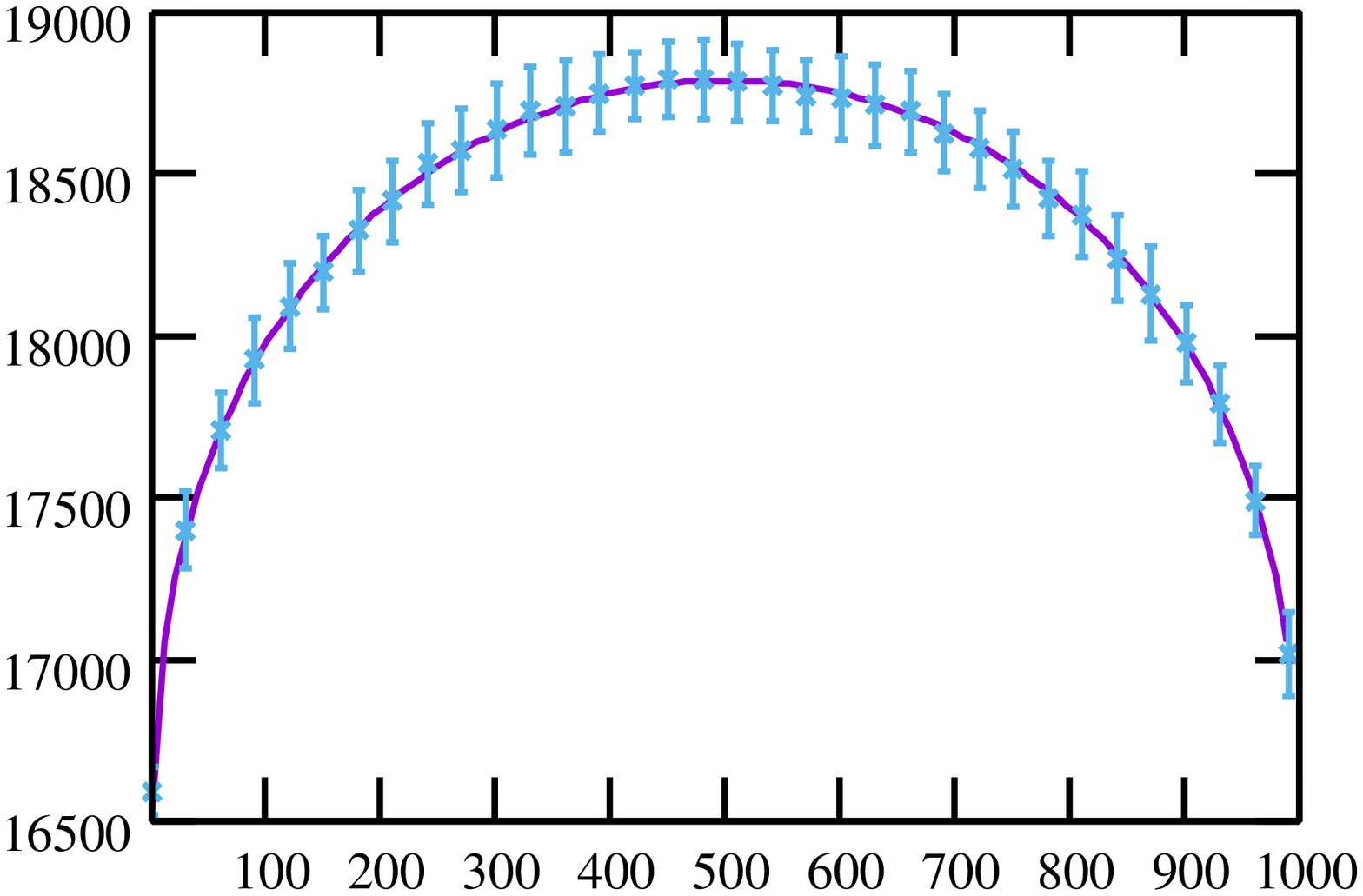} 
\caption{\label{fig-chi-sq-600}$M=1000$ and $N=600$, chi-squared distribution.{The horizontal axis shows the position of the fiducial machine 
 $\mu$ and the vertical axis shows the makespan $\makespan(\mu)$.}
}
\vspace{0.5cm}
\includegraphics[width=0.8\hsize]{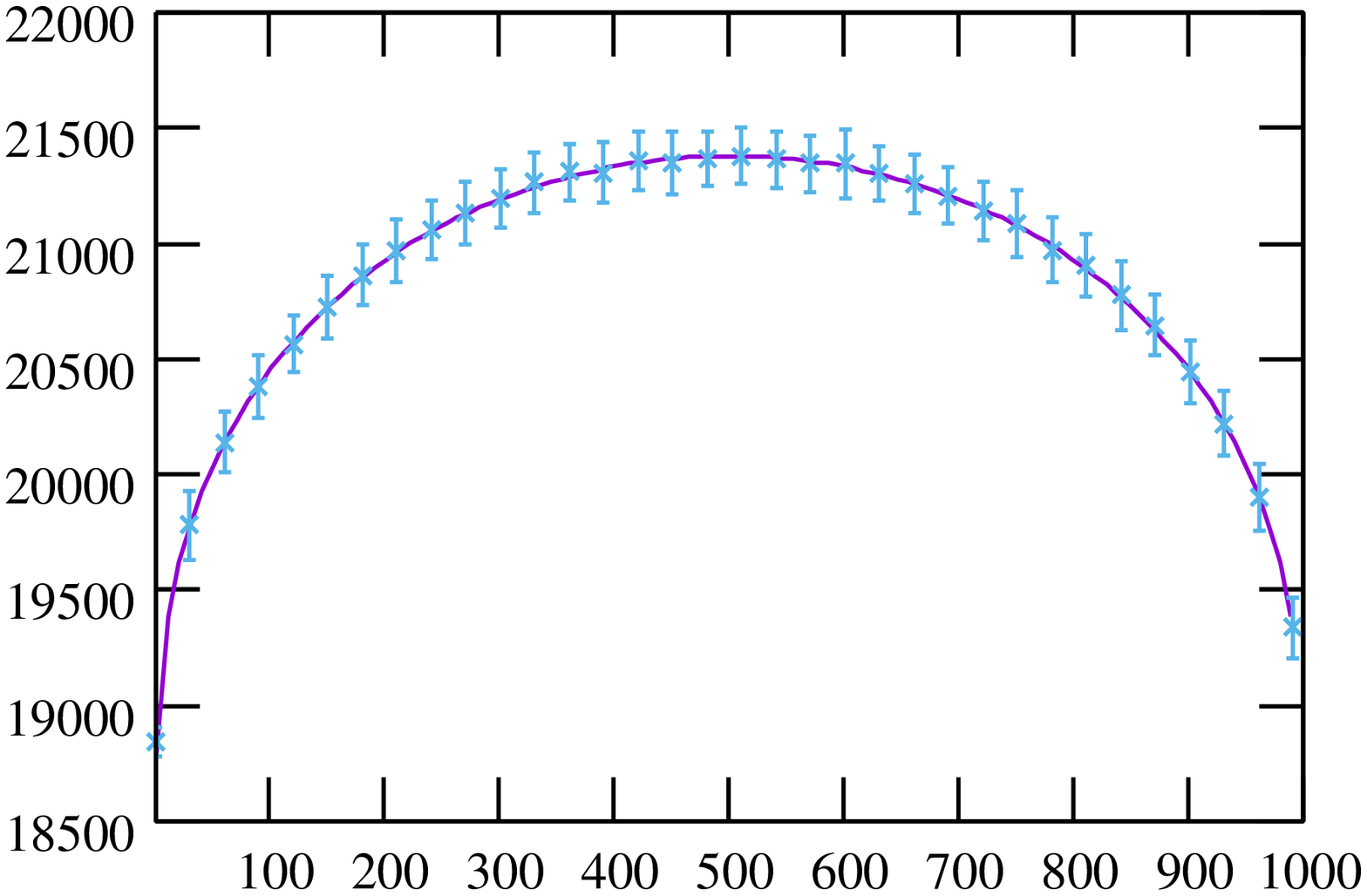} 
\caption{\label{fig-chi-sq-800}$M=1000$ and $N=800$, chi-squared distribution.{The horizontal axis shows the position of the fiducial machine 
 $\mu$ and the vertical axis shows the makespan $\makespan(\mu)$.}
}

\end{center}
\end{figure}

Lastly, the relationship between the average of makespan and the 
position of the fiducial machine in the {processing system} is examined for the case that 
the processing time is independently and identically distributed according to 
the chi-squared distribution
\bea
f(X_{\mu,i})\eq
\left\{
\begin{array}{ll}
\f{(X_{\mu,i})^{\f{k}{2}-1}}{2^{\f{k}{2}}\Gamma(k/2)}e^{-\f{X_{\mu,i}}{2}}
&X_{\mu,i}>0\\
0&\text{otherwise}
\end{array}
\right.,
\eea
where ${\rm E}_X[X_{\mu,i}]=k$ and 
${\rm Var}[X_{\mu,i}]={\rm E}_X[X^2_{\mu,i}]-({\rm 
E}_X[X_{\mu,i}])^2=2k$. If $k=7$, then 
the mean and variance are consistent with {those for the previous uniform distributions}. 
Fig. \ref{fig-chi-sq-200} to Fig. \ref{fig-chi-sq-800} show that the expectation of makespan is concave 
for this case and similar to the other distributions considered up to a 
statistical fluctuation. Namely, from these numerical experiments, we can assume that 
the asymptotical behavior of makespan is 
universal, that is, it does not depend on the distribution of processing time.
The next section will discuss this universality  mathematically.

\begin{table}[t]
\begin{center}
\begin{tabular}{|c|c|c|c|c|l|}
\hline
$M$&$N$&$A$&$B$&$\a$&Figure\\
\hline
$1000$&$200$
&$89.9326$
&$8289.92$
&$0.51588$
&Fig.\ref{fig-chi-sq-200}
\\
\hline
$1000$&$400$
&$150.0050$
&$9884.91$
&$0.48407$
&Fig.\ref{fig-chi-sq-400}
\\
\hline
$1000$&$600$
&$177.9403$
&$11131.28$
&$0.49390$
&Fig.\ref{fig-chi-sq-600}
\\
\hline
$1000$&$800$
&$211.6288$
&$12650.38$
&$0.48702$
&Fig.\ref{fig-chi-sq-800}
\\
\hline
\end{tabular}
\caption{Parameter estimates for fitting functions (plotted in the figures) for the case that the processing time is chi-squared distributed.\label{tab-chi-sq}}
\end{center}
\end{table}
\clearpage

\section{Mathematical Discussion\label{mathematical-discussion}
\label{sec4}}
In this section, we discuss the function $h$ appearing in \Sref{eq17} or its 
asymptotic function $g$, defined in detail below, for large numbers of 
machines and jobs.
Unexpectedly, $g$ is related to a function called the {\it shape function} in percolation theory, about which much has already been revealed. We will see how 
{percolation theorems can be used to derive
scheduling theorems}.

\subsection{Percolation Theory and Shape Function}
Let discuss the forward scheduling \Sref{eq5} and \Sref{eq6}. The {following argument for forward scheduling}  can be applied to backward scheduling, by symmetry.
Set $T_{0,0}=0$. Then for $(m,n) \in \mathbb{Z}_{\geq 0}^2$, we can show
\bea
T_{m,n}\eq\max_{\gamma\in \{(0,0)\rightsquigarrow(m,n)\}} \sum_{(\mu,i)\in\gamma}  X_{\mu,i}.
\label{eq:s}
\eea
Here, $\{(0,0)\rightsquigarrow(m,n)\}$ means the set of the shortest paths in the square grid graph connecting $(0,0)$ and $(m,n)$.
This relation can be derived from \Sref{eq5} and \Sref{eq6} by induction.

If $X_{\mu,i}\in \mathbb{R}$ for $\mu,i\in \mathbb{Z}_{\geq 0}$ is sampled i.i.d. from a probability distribution $F$, then the model above is known as the {\it two-dimensional 
last-passage directed site percolation}. This name comes from the problem of finding the time of 
last-passage {of a particle which starts from the origin and moves} only in the ``right'' or ``up'' direction for 
a lattice with weights {(time required to pass through)} on its vertices  (rather than its edges). {Last}-passage  directed site percolation has been studied mainly in mathematics and physics 
\cite{corwin-2012,johansson-2000,martin-2004}, and it is also referred to as {\it zero-temperature directed polymer in a random environment} when used as a polymer model.

Asymptotic properties of $T_{m,n}$ have been reported (see 
\cite{martin-2004}). Throughout the present paper, we assume ${\rm 
E}_X[X_{\mu,i}]<\infty$ and that $X_{\mu,i}$ is non-degenerate, i.e., ${\rm Var}[X_{\mu,i}]>0$ for any $\mu$ and $i$. There exists a function $g(\vec{v})$ on $\mathbb{R}^2_{\geq 0}$ as an almost sure asymptotic limit of $T$: for every $\vec{v}=(v_1,v_2)\in \mathbb{R}^2_{\geq 0}$,
\bea
\label{eq35}
\frac{1}{N} T_{\lfloor N v_1\rfloor,\lfloor N v_2\rfloor} ~\mathop{\rightarrow}^{\rm a.s.}~ g(v_1,v_2)
\label{eq:t2g}
\label{eq32}
\eea
as $N\rightarrow \infty$. 
The convergence in $L_1$ also holds:
\bea
\label{eq33}
\lim_{N\rightarrow \infty} \frac{1}{N} {\rm E}_X[T_{\lfloor N v_1\rfloor,\lfloor N v_2\rfloor}] = g(v_1,v_2).
\label{eq:Et2g}
\eea

Here, $\displaystyle\mathop{\rightarrow}^{\rm a.s.}$ means almost sure convergence
(in other words, convergence with probability one)
and $\lfloor \cdot \rfloor$ is the floor function.
The function $g(v_1,v_2)$ is called a {\it shape function} or is sometimes known as a {\it time constant}.
Note that the existence of such a function $g$ is not obvious but can be proved using Kingman's subadditive ergodic theorem (see \cite{martin-2004}). 

\subsection{Shape Functions and Makespan of Scheduling}
Let $\lambda$ be the mean of random variable $X_{\mu,i}$ and $F(X_{\mu,i})$ be its cumulative 
distribution function. 
To make the argument simpler, we will use the normalization {$\tilde{X}_{\mu,i} = 
X_{\mu,i}-\lambda$.
The corresponding $T_{m,n}$ becomes}
\bea
\tilde{T}_{m,n} \eq \max_{\gamma\in \{(0,0)\rightsquigarrow(m,n)\}} \sum_{(\mu,i)\in\gamma}  (x_{\mu,i}-\lambda) \nn
\eq T_{m,n}-(m+n)\lambda.
\eea
Therefore, the corresponding shape function $\tilde{g}$ is represented as 
\bea
\tilde{g}(v_1,v_2) \eq g(v_1,v_2)-\lambda (v_1+v_2).\label{eq:g_tilde}
\eea

Next we derive  makespan for the algorithm in Section \ref{sec2.5}.  Assume 
$M$, the number of machines, and $\nu$, the index of the fiducial 
machine, 
are proportional to  $N$ in the following sense:
\bea
M\eq\kappa N + o(N), \label{eq:M_N}
\eea
and
\bea
\nu\eq\tau M + o(M). \label{eq:nu_M}
\eea
for some $\kappa\geq 0$ and $\tau\in [0,1]$. 
 The following theorem states 
that the normalized makespan is asymptotically represented by the 
sum of two shape functions.

\begin{thm}
\label{thm:makespan}
Assume \Sref{eq:M_N} and \Sref{eq:nu_M}. Then
\bea
\label{eq38}
\frac{1}{N} T(\nu) 
&\mathop{\rightarrow}^{\rm a.s.} g(\kappa(1-\tau),1)+g(\kappa\tau,1) -\lambda \label{eq:thm:makespan1}\\
\label{eq39}
&= \tilde{g}(\kappa(1-\tau),1)+\tilde{g}(\kappa\tau,1) + \lambda(\kappa+1) \label{eq:thm:makespan2}
\eea
and
\bea
\label{eq40}
\frac{1}{N} {\rm E}_X[T(\nu)]
&\rightarrow g(\kappa(1-\tau),1)+g(\kappa\tau,1) -\lambda \label{eq:thm:makespan3}\\
\label{eq41}
&= \tilde{g}(\kappa(1-\tau),1)+\tilde{g}(\kappa\tau,1) + \lambda(\kappa+1) \label{eq:thm:makespan4}
\eea
as $N\rightarrow \infty$.
\end{thm}
The proof is in \ref{appendix:A}.

From this theorem, the asymptotic values of the normalized makespan are obtained as
\bea
h_{\kappa}(\tau):=g(\kappa(1-\tau),1)+g(\kappa\tau,1) -\lambda
\eea
and we can prove the following property of $h_{\kappa}$.

\begin{thm}
\label{thm:4.2}
Assume the function $h_{\kappa}$ on $[0,1]$ is a concave function which attains its maximum at the midpoint $\tau=1/2$ and
its minimum at the endpoints $\tau=0$ and $1$.
\end{thm}
The proof is in \ref{appendix:B}.

\subsection{Makespan for the Exponential or Geometric Distribution}
For only a few types of distributions $F$ are the shape functions $g$ known.
The following theorem is {based on} the results proved in \cite{johansson-2000}.

\begin{thm} \label{thm:sqrt}
If $F$ is the cumulative distribution function of an exponential distribution or a geometric distribution with mean $\lambda$
and standard deviation $\sigma$, then
\bea
\label{eq37}
g(\xi,1)=\lambda(1+\xi) + 2\sigma\sqrt{\xi}  \label{eq:shape-sqrt}
\eea
and
\bea
\label{eq38}
\tilde{g}(\xi,1)=2\sigma\sqrt{\xi}. \label{eq:shape-sqrt-tilde}
\eea
\end{thm}
In \ref{appendix:C}, we note which results in \cite{johansson-2000}
correspond to Theorem \ref{thm:sqrt}.

From Theorem \ref{thm:makespan} and Theorem \ref{thm:sqrt}, the normalized makespan 
for the exponential distribution and the geometric distribution is obtained as follows.
\begin{cor}
Assume \Sref{eq:M_N} and \Sref{eq:nu_M}. 
If $F$ is the exponential distribution or geometric distribution with mean $\lambda$
and standard deviation $\sigma$, then
\begin{align}
\frac{1}{N} T(\nu) 
&\mathop{\rightarrow}^{\rm a.s.} 
2\sigma\left(\sqrt{\kappa(1-\tau)}+\sqrt{\kappa\tau}\right) + \lambda(\kappa+1)
\end{align}
and
\begin{align}
\frac{1}{N} {\rm E}_X[T(\nu)] \rightarrow
2\sigma\left(\sqrt{\kappa(1-\tau)}+\sqrt{\kappa\tau}\right) + \lambda(\kappa+1)
\end{align}
as $N\rightarrow \infty$.
\end{cor}

In \cite{martin-2004}, they noted that the exponential or geometric distribution $F$ ``is essentially the only
nontrivial case (whether directed or undirected, first- or last-passage) where
the form of the shape function $g$ above is known.''
As far as we know, the corresponding makespan function is also unknown except for the case of an exponential or
geometric distribution.
Remark that, on the other hand, they proved that $g(\epsilon,1)$ behaves as $\sigma\sqrt{\epsilon} + \lambda(1+\epsilon)$ for small $\epsilon>0$ for a general distribution $F$.
Therefore, the normalized makespan function $\frac{1}{N}T(\nu)$ behaves asymptotically 
parabolic for sufficiently small $\nu/N$.

\section{Conclusions and Future Work}
In this paper, 
we have discussed the relationship between the asymptotical behavior 
of makespan and the position of the fiducial machine in a {processing system} 
both numerically and mathematically. When 
a job processing order is fixed, 
the first and last machines in this {processing system} 
minimize the makespan 
and the center machine in the {processing system} maximizes the makespan instead of the 
bottleneck machine. Regarding  the dispatching rules discussed in 
several works, shortest processing time and 
longest processing time, 
the makespan of 
the job processing order without a dispatching rule 
is hardly distinguishable from the makespan of 
the job processing order based on the shortest processing time rule, 
while it is less than the makespan of 
the job processing order based on the longest processing time rule.
Namely, the dispatching rules do not work well for the flowshop scheduling {
problem. Also,} the makespan is strongly influenced by the position of 
the fiducial machine in the {processing system}. Moreover, 
numerical findings are supported by  the property  of shape functions discussed in 
percolation theory. In addition, 
when the processing time is independently and identically distributed 
according to an exponential or geometric distribution, it is well known that the shape function has  
a term proportional to the square of the position of the fiducial machine, although 
the shape functions in the cases that the processing time is independently and identically distributed according to a discrete uniform distribution, continuous uniform distribution, or
chi-squared distribution 
might have a term proportional to the square of the position, consistent with the numerical results of the \shape function.

In the future work, because our results here regarding the shortest processing time and longest processing time 
as the dispatching rules show that not all dispatching rules consistently work well 
with the flowshop scheduling problem, we need to compare 
other dispatching rules, both numerically and mathematically. 
Next, 
we considered here only i.i.d. processing times in order to simplify our discussion;
however, since 
processing times in practice are weakly or strongly correlated with 
each other, 
the makespan in
the case of 
correlated processing times needs to be examined in detail.

\section*{Acknowledgement}
One of the authors (TS) appreciates the fruitful comments of H. Yamamoto, 
I. Arizono, 
and Y. Kajihara. 
The work is partly supported by Grants-in-Aid Nos. 24700288, 26280009, 
24710169, and 15K20999; 
JST PRESTO; the
President Project for Young Scientists at Akita Prefectural University; 
research project No. 50 of the National Institute of Informatics, Japan; 
research project No. 5 of the Japan Institute of Life Insurance; 
research project of the Institute of Economic Research Foundation; Kyoto University; 
research project No. 1414 of the Zengin Foundation for Studies in Economics 
and Finance; 
research project No. 2068 of the Institute of Statistical Mathematics; 
research project of Mitsubishi UFJ Trust Scholarship Foundation; 
and research project No. 2 of the Kampo Foundation.

\appendix
\section{Proof of Theorem \ref{thm:makespan}}
\label{appendix:A}

Any shape function $g$ satisfies the following properties (listed in \cite{martin-2004}):
\begin{description}
\item[(P1)] $g(0,0)=0$,
\item[(P2)] $g(v_1,v_2)=g(v_2,v_1)$ for $v_1,v_2\in \mathbb{R}_{\geq 0}$,
\item[(P3)] $g(\alpha \vec{v})=\alpha g(\vec{v})$ for $\alpha\geq 0$,
$\vec{v}\in \mathbb{R}^2_{\geq 0}$, 
\item[(P4)] $g(\vec{v})+g(\vec{w})\geq g(\vec{v}+\vec{w})$ for $\vec{v},\vec{w}\in \mathbb{R}^2_{\geq 0}$.
\end{description}
By \Sref{eq35},
$\tilde{g}(v_1,v_2)$ also satisfies properties (P1)-(P4).
Note that by the law of large numbers,
\bea
\frac{1}{N} T_{\lfloor N v\rfloor,0}\eq
\frac{1}{N}\sum_{\mu=0}^{\lfloor N v\rfloor} X_{\mu,0}
\mathop{\rightarrow}^{\rm a.s.}~ v\lambda.
\eea
By comparing this and \Sref{eq:t2g}, we obtain the properties
\begin{description}
\item[(P5)] $g(v,0)=g(0,v)=v\lambda$ for $v\geq 0$ and
\item[(P5)'] $\tilde{g}(v,0)=\tilde{g}(0,v)=0$ for $v\geq 0$.
\end{description}

Normalizing the makespan $T(\nu)$ by the number of jobs $N$, 
\bea
\f{1}{N} T(\nu)\eq
\f{1}{N}\{t_{M,N}(\nu)-s_{1,1}(\nu)\} \nonumber\\
\eq\f{1}{N}\{t_{M,N}(\nu)-s_{\nu,1}(\nu)\} - \f{1}{N} \{t_{\nu,N}(\nu)-s_{\nu,1}(\nu)\}\nonumber\\
&&+\f{1}{N}\{t_{\nu,N}(\nu)-s_{1,1}(\nu)\}. \label{eq:T}
\eea
The first and the third terms follow the same distributions as $T_{M-\nu,N-1}$ and $T_{\nu-1,N-1}$, respectively,
where $T_{M-\nu,N-1}$ and $T_{\nu-1,N-1}$ are two independent random variables generated by \Sref{eq:s}.

Then by \Sref{eq:t2g},
\bea
\frac{1}{N}(T_{M-\nu,N-1} + T_{\nu-1,N-1}) ~\mathop{\rightarrow}^{\rm a.s.}~
g(\kappa(1-\tau),1)+g(\kappa\tau,1).
\eea

The second term of the right-hand side of \Sref{eq:T} converges to $-\lambda$ since
\bea
\frac{1}{N} \{t_{\nu,N}(\nu)-s_{\nu,1}(\nu)\}\eq
\frac{1}{N} \sum_{i=1}^N X_{\nu,i}
\mathop{\rightarrow}^{\rm a.s.} \lambda
\eea
by the strong law of large numbers.
Thus, \Sref{eq38} and \Sref{eq39} of Theorem \ref{thm:makespan} hold.
Similarly but by using \Sref{eq33} instead of \Sref{eq32},
\Sref{eq40} and \Sref{eq41} can be proved.

\section{Proof of Theorem \ref{thm:4.2}}
\label{appendix:B}
It is sufficient to prove that $\tilde{g}(\xi,1)$ is a strictly increasing concave function of $\xi\geq 0$.
By properties (P3) and (P4), $g(\xi,1)$ is concave and so is $\tilde{g}(\xi,1)$.
Because of (P5)' and the concavity, $\tilde{g}(\vec{v})>0$ for 
$\vec{v}\in \mathbb{R}_{> 0}^2$ or $\tilde{g}(\vec{v})\equiv 0$ for $\vec{v}\in \mathbb{R}_{\geq 0}^2$.
Next assume that  $X_{\mu,i}$ are non-degenerate random variables and 
${\rm E}_X[X_{\mu,i}]=0$. Then,\bea
{\rm E}_X[T_{1,1}]\eq{\rm E}_X[\max(X_{1,2},X_{2,1})]\nn
\eq{\rm E}_X[X_{2,1}+\max(0,X_{1,2}-X_{2,1})]\nn
&>&0.
\eea
Since ${\rm E}_X[T_{N,N}]\geq N {\rm E}_X[T_{1,1}]$, 
\bea
\tilde{g}(1,1)=\lim_{N\rightarrow \infty} \frac{1}{N}{\rm E}_X[T_{N,N}]\geq {\rm E}_X[T_{1,1}]>0.
\eea
Thus, $\tilde{g}(\vec{v}) \not\equiv 0$, which implies $\tilde{g}(\vec{v})>0$ for 
$\vec{v}\in \mathbb{R}_{> 0}^2$.

For $0\leq \xi<1$ and $0<\epsilon\leq 1-\xi$, the point
\bea
\vec{p}=\left(\frac{1+\xi}{1+\xi+\epsilon}(\xi+\epsilon),\frac{1+\xi}{1+\xi+\epsilon}\right)
\eea
is an {internal} division of {the line segment joining points} $(\xi,1)$ and $(1,\xi)$. Thus by the concavity and
and property (P3) of $\tilde{g}$,
\begin{align}
\frac{1+\xi}{1+\xi+\epsilon} \tilde{g}(\xi+\epsilon,1)
&=\tilde{g}\left(\frac{1+\xi}{1+\xi+\epsilon}(\xi+\epsilon),\frac{1+\xi}{1+\xi+\epsilon}\right) \nonumber\\
&\geq \tilde{g}(\xi,1)~ (= \tilde{g}(1,\xi)).
\end{align}
Therefore, $\tilde{g}(\xi+\epsilon,1)>\tilde{g}(\xi,1)$ and $\tilde{g}(\xi,1)$ is a strictly increasing function of $\xi\geq 0$.  
\qed

Note that $g$ is continuous on $\mathbb{R}^2_{> 0}$ since it is concave, and, therefore, $h_\kappa$ is continuous for $\tau>0$.
In \cite{martin-2004}, they proved that $g$ is continuous on $\mathbb{R}^2_{\geq 0}$, the whole domain including the boundary, if distribution $F$ satisfies
\bea
\int^\infty_0 (1-F(s))^{1/2} \diff s < \infty.
\label{eq:F1}
\eea
Remark that condition (\Sref{eq:F1}) is satisfied for almost all common distributions
since it corresponds to the density functions decaying faster than in cubic order.

\section{Derivation of Theorem \ref{thm:sqrt}}
\label{appendix:C}
Theorem \ref{thm:sqrt} is included in the results proved in 
\cite{johansson-2000}. Here, we will note which results in 
\cite{johansson-2000} correspond to Theorem \ref{thm:sqrt}. 

From Theorem 1.1 of \cite{johansson-2000}, the geometric distribution
$P(X_{\mu,i}=k)=(1-q)^k q$ for $k\in \mathbb{Z}_{\geq 0}$,
whose mean is ${\rm E}_X[X_{\mu,i}]= \frac{q}{1-q}$
and variance is $V[X_{\mu,i}]= \frac{q}{(1-q)^2}$,
induces the shape function
\bea
g(\xi,1)=\frac{(1+\sqrt{q\xi})^2}{1-q}-1=\frac{q}{1-q}(1+\xi)+2\sqrt{\frac{q}{(1-q)^2}\xi}.
\eea
This coincides with \Sref{eq:shape-sqrt}.

From Theorem 1.6 of \cite{johansson-2000}, the exponential distribution whose mean is one
(and, therefore, variance is one) induces the shape function
\bea
g(\xi,1)=\left(1+\sqrt{\xi}\right)^2=1+\xi+2\sqrt{\xi}.
\eea
Since the exponential distribution with parameter $\lambda$,
whose mean is $\lambda$ and variance is $\lambda^2$, is obtained by multiplying $X_{\mu,i}$ by 
$\lambda$ and ,
the shape function becomes
\bea
g(\xi,1)=\lambda(1+\xi)+2\lambda\sqrt{\xi}.
\eea
This also coincides with \Sref{eq:shape-sqrt}.

Remark that in Theorem 1.1 and Theorem 1.6 of \cite{johansson-2000}, $\xi\geq 1$ is assumed 
but {statements of these theorems} can be generalized to $\xi>0$ by setting $\xi:=1/\xi$. The theorem is 
easy to prove for $\xi=0$.
Eq.(\ref{eq38}) for $\tilde{g}$ is a direct consequence of
\Sref{eq37}.

\bibliographystyle{amsplain}

\end{document}